\documentclass[sn-mathphys-num]{sn-jnl}

\usepackage{lscape}
\usepackage{longtable}
\usepackage{graphicx}%
\usepackage{multirow}%
\usepackage{amsmath,amssymb,amsfonts}%
\usepackage{amsthm}%
\usepackage{mathrsfs}%
\usepackage[title]{appendix}%
\usepackage{xcolor}%
\usepackage{textcomp}%
\usepackage{manyfoot}%
\usepackage{booktabs}%
\usepackage{algorithm}%
\usepackage{algorithmicx}%
\usepackage{algpseudocode}%
\usepackage{listings}%
\usepackage{longtable}

\theoremstyle{thmstyleone}%
\newtheorem{theorem}{Theorem}
\theoremstyle{thmstyletwo}%
\newtheorem{remark}{Remark}%

\theoremstyle{thmstylethree}%
\newtheorem{definition}{Definition}%
\newtheorem{lemma}{Lemma}%
\begin{document}

\title[Article Title]{Bi-objective Two Phase  Transportation Problem}

 
\author[1]{\fnm{Prabhjot } \sur{Kaur}}\email{prabhjotkaur@pu.ac.in}

\author*[1]{\fnm{ Kalpana } \sur{Dahiya}}\email{kalpanas@pu.ac.in}
 
\author[3]{\fnm{ Vanita } \sur{Verma}}\email{v\_verma1@yahoo.com}
 
\affil[1]{\orgdiv{University Institute of Engineering and Technology}, \orgname{Panjab University}, \orgaddress{  \city{Chandigarh}, \postcode{160014}, \state{Chandigarh}, \country{India}}}

\affil[3]{\orgdiv{Department of Mathematics}, \orgname{Panjab University}, \orgaddress{  \city{Chandigarh}, \postcode{160014}, \state{Chandigarh}, \country{India}}}


\abstract{This paper discusses a bi-objective two-phase transportation problem (BTPTP) in which transportation takes place in two phases.  The whole  set of source-destination links is partitioned into two disjoint sets namely, Phase-I and Phase-II. The transportation among the source-destination links of Phase-II set is done only after the transportation among the source-destination links of Phase-I set is complete. The transportation among the source-destination links in each phase  is done in parallel. The problem discussed in this paper  concentrates on minimizing two objectives viz., sum of transportation times and sum of transportation costs of both the phases, simultaneously.
	Due to conflicting nature of   the objectives, a polynomial time iterative algorithm named as BTPTP-Algorithm is developed that finds all of its non-dominated solutions.  At each iteration, the BTPTP-Algorithm solves some specific  restricted cost minimizing transportation problems and finds  corresponding optimal feasible solutions. The restrictions imposed on the source-destination links of Phase-I and Phase-II sets, in each of these problems, depend upon their corresponding Phase-I and Phase-II transportation times obtained at the previous step. The algorithm is developed in such a way that it expels all the dominated solutions and records all the non-dominated  solutions  of the problem, systematically. A proof of the BTPTP-Algorithm's capability of recording all the  non-dominated   points is provided. Further, a numerical illustration is given in support of the theory and  computational behavior of the BTPTP-Algorithm for various BTPTP instances of different sizes, is provided in terms of CPU time.}

\keywords{Transportation, Multi-objective, Pareto optimal  solutions, Non-dominated points, Two-phase }

\pacs[MSC Classification]{90C26   90C27}

\maketitle
 
\section*{Funding} The work presented in this paper is supported by grant from the  Science and Engineering Research Board, Government of India (File no. MTR/2019/000723).
\section*{Conflict of Interest/Competing Interest} The authors  declare that they have no conflict of interest.
\section*{Availability of data and material} Not applicable
\section*{Code availabilty} Not applicable	
\section*{JEL Codes} C61, C63 \section*{Acknowledgements} Not applicable
\section*{List of Abbreviations}
The list of abbreviations for the terms frequently used in the following  text is given below:
$$\begin{array}{lll}
	\mbox{BTPTP}&:&\mbox{bi-objective two-phase transportation problem}\\
	\mbox{OFS}&:&\mbox{Optimal Feasible Solution }\\
	\mbox{CMTP}&:&\mbox{Cost Minimizing Transportation Problem }\\
	\mbox{OMFS}&:&\mbox{Optimal M-feasible Solution}\\
\end{array}$$
	\section{Introduction}
In many real world problems, it becomes essential to optimize more than one objectives simultaneously.   For example, while purchasing a car, one would make a decision according to the price, comfort, fuel consumption and power of each of the alternatives. Obviously,  a powerful car which is cheap, comfortable  and has low fuel consumption, is preferred.  Ideally, it is unusual to expect that a cheapest car is the most comfortable, economic  and powerful one. That means the best choice have to be made out of   some qualitative alternatives available. These kind of optimization problems are called multi-criteria optimization problems.
In such problems, there are generally, no decision variables (only alternatives are available) due to which these problems cannot be formulated, mathematically. However, if in an optimization problem with multiple objectives, each objective  can be formulated  mathematically, in terms of some decision variables and  has to be   minimized or maximized,  such problems belong to the class of multi-objective optimization problems.   Multi-objective optimization problems find their applications in many fields of science, engineering and economics.
In a trivial case of multi-objective optimization problem, there always exists a point in the solution set that simultaneously optimizes all the objective functions. However,
if the multiple objectives of a multi-objective optimization problem  are of conflicting  nature, they may not be optimized simultaneously. In this case, one has to search for a compromising or a trade-off solution that involves loss of one objective value in return for the gains in the others. Such solutions  form the so called Pareto front and  called Pareto optimal solutions. A Pareto optimal solution gives rise to a non-dominated point of the problem. Researchers study multi-objective optimization problems from different viewpoints and developed different solution methodologies to find their Pareto optimal solutions.  Some of the recent work done on multi-objective optimization problems in which various  numerical and   evolutionary algorithms are developed using crisp and fuzzy  data can be found in \cite{Bazine,Branke,Joydeep,Fakhar,Ghaznavi,bal1,bal2,Tan,bal}.
Ultimately, the decision to quantificate all the Pareto optimal solutions of a multi-objective optimization problem along the Pareto front imposes considerable cognitive burden on the human decision-makers who have to choose one solution from a large number of Pareto optimal solutions. Their choice of this particular solution depends on the additional subjective preference information such as the weight associated with each objective (if any) or the interval specified for the value of an objective function to lie within. Without such additional subjective preference information, all the Pareto optimal solutions are considered equally good.
\par In real world engineering problems also, many challenges deal with multiple objectives instead of a single objective. Accordingly, the study of transportation problems with more than one objective becomes important. In literature, multi-objective transportation problems have been intensively investigated by several researchers  \cite{Climaco,Current,Diaz,Diaz2,Isermann,Kasana,Lee,Nomani,Ringuest} and various solution strategies have been developed to solve them. A comparison of the objectives considered by them   is given in  Table \ref{1gg} by writting a `yes' for the  active objective  in a particular  problem. However, a brief explanation of some  of the important papers is given as follows:
\par Lee and Moore \cite{Lee} discussed the optimization of transportation problems with multiple objectives such as fulfillment of transportation schedule contracts and minimizing transportation hazards using goal programming approach.
Diaz \cite{Diaz,Diaz2} presented an alternative procedure to generate all non-dominated solutions to the  linear  multi-objective transportation problem followed by a complete description of all of its Pareto optimal solutions. This approach depends upon specifying an \textit{`a priori'}  measure of the closeness of any compromise solution to the ideal solution.   Current et al. \cite{Current} also gave a review of multi-objective design of transportation networks in 1986. In 1987,  Ringuest and Rinks \cite{Ringuest}  presented  interactive solutions for linear multi-objective transportation problems. Then, Bit et al. \cite{Bit1} presented an application of fuzzy linear programming to the linear multi-objective transportation problem to obtain  efficient solutions as well as an optimal compromise solution.
They  \cite{Bit1}  suggested that the  interactive algorithms developed by  Diaz \cite{Diaz} and Ringuest and Rinks  \cite{Ringuest}   are only applicable to particular types of the multi-objective transportation problems, however,  the fuzzy programming algorithm is applicable to all types of multi-objective transportation problems.  Later on, Bit et al. \cite{Bit} also gave a fuzzy programming approach to multi-objective solid transportation problem.    Over the last three decades,  multi-objective transportation problems have been widely studied under usual as well as uncertain environment  \cite{Waiel,Abd,Climaco,Ehrgott,Kasana,Lushu,Nomani,Verma,Yu}.	    	 
\par Further, various transportation problems with two objectives (such as transportation time, cost, deterioration, efficiency etc., taken two at a time)  have also attained due attention of various researchers \cite{103,43,105,102}. Some  variants of transportation problems such as fixed charge transportation problem, bulk transportation problems have also been discussed thoroughly in literature (See Table  \ref{1gg}). A brief review of these bi-criteria transportation problems is given as follows:
\par
Time-cost trade-off in fixed-charge bi-objective  transportation problem have been discussed by Basu et al. \cite{Basu} whereas Gupta et al. \cite{108} discussed  time-cost trade-off relations in a bulk transportation problem. Further, keeping in view the realistic market scenario which advocates that fact that a transportation problem can not always be studied with crisp and  certain data,  Mitsuo et al. \cite{40} developed a genetic algorithm to discuss a bi-objective  solid transportation problem  with fuzzy numbers in 1997.  An algorithm for solving a fixed charge bi-objective transportation problem with restricted flow was given by Thirvan et al. \cite{Thirwan}.
Pandian and Anuradha  \cite{Pandian} discussed another bi-objective transportation problem in which the percentage level of satisfaction for a solution of the transportation problem is introduced. Furthermore,  a method to find non-dominated solutions for bi-objective integer transportation problem was given by Basirzadeh and Delfzergan \cite{Basirzadeh}.
\begin{landscape}
	\begin{table}\centering
		\caption{Literature review of Multi-objective/Bi-objective transportation problems}
		\begin{tabular}{lcccccccc}\label{1gg}
			& \textbf{ Fullfilling } &  & \textbf{Multiple } & & & &&\\
			\textbf{Author/objectives} $\rightarrow$    & \textbf{ transportation  } & \textbf{hazards/} & \textbf{objectives}&\textbf{fuzzy}& \textbf{Cost }& \textbf{interval}&\textbf{Time}&\textbf{Efficiency}\\
			$\downarrow~~~~~~~~~~~~~~~~~~~~~~$  & \textbf{  schedule} & \textbf{deterioration}&\textbf{  }&  & &\textbf{parameters}& &\\
			& \textbf{ contract  } &&\textbf{  }  &&&&& \\
			\hline
			Lee and Moore \cite{Lee}  &  yes    &yes  &x&x&x&x&x&x\\
			\hline
			Isermann \cite{Isermann}&x  &x   & yes &x&x&x&x&x\\
			\hline
			Diaz \cite{Diaz,Diaz2} &x&x&yes& x  &x&x&x&x\\\hline
			Current et al. \cite{Current}&x&x&yes&x   &x&x&x&x\\\hline
			Ringuest and Rinks \cite{Ringuest}  &x&x&yes&   x&x&x&x&x\\\hline
			Bit et al. \cite{Bit1},Nomani et al. \cite{Nomani}&x&x&yes&  yes &x&x&x&x\\\hline
			Chanas et al. \cite{Chanas}&x&x&x&yes&x&x&x\\
			Chanas and  Kuchta \cite{Chanas1},  &x&x&x&yes&yes&x&x&x\\\hline
			Mitsuo et al. \cite{40,104}&x&x&x&yes&yes&x&x&x\\\hline
			Verma et al. \cite{Verma} &x&x&yes&yes& x&&x&x\\\hline
			Waiel et al. \cite{Waiel} &x&x&yes&yes&x&x&x&x \\\hline
			
			Yu et al. \cite{Yu}     &x&x&yes& & yes&yes&x&x\\\hline

			Gupta et al. \cite{105}, Prasad et al. \cite{42}&x&x&x & x  &yes & x&yes&x\\
			Quddoos et al. \cite{Quddoos2,Quddoos},     Singh et al. \cite{Singh2018}&x&x&x & x  &yes &x &yes&x\\
			Basu et al.  \cite{Basu}, Gupta et al. \cite{108}  &x&x&x &x   &yes &x &yes&x\\			 
			Aneja and Nair \cite{103}&x&yes&x & x  &yes &x &x&x \\\hline
			Malhotra \cite{102} &x&yes&x & x  &yes &x &yes &x\\\hline
			Derigs \cite{43}&x&x & x&  x &yes & x&yes &yes\\\hline
	\end{tabular} \end{table}
\end{landscape}
\par It is usually assumed that in a transportation problem, the demands of the destinations are satisfied by transporting the goods from various sources in one go. However, in many real life situations, due to some budgetary or storage constraints, it is not always possible and therefore, the transportation has to be done in many stages thus, giving  rise to   multi-stage transportation problems.
Table \ref{tt} provides a brief review of the work done on various two-stage/multi-stage transportation problems.  A very important study on sustainable closed-loop supply chain with efficiency and resilience systematic framework by developing a mathematical model using two-stage, mixed-integer linear programming has been done by Mehrjerdi and Lotfi \cite{Mehrjerdi}.
\par In this paper,    a two-phase transportation problem is discussed. It is a specific case of the two-stage transportation problems in which, the whole set of source-destination links is divided into two mutually disjoint and exhaustive sets, viz., Phase-I and Phase-II sets.  In Phase-I of the problem, the transportation is carried out only along the links belonging to Phase-I set and after the completion of   Phase-I transportation, only the links belonging to  Phase-II set are used for transportation.
The difference between a two-stage  and a two-phase transportation  problem is that in the former problem, the transportation is allowed on all the source-destination links in both the stages  whereas in the latter problem, the source-destination
links that are allowed to participate in transportation process in Phase-I of the problem are blocked in the transportation process of Phase-II and vice-versa. Some of the two-stage transportation problems have also been studied in literature  as the problems of distributing   a homogeneous product  from   plants to distribution centers and then to customers \cite{Calvete,Calvete1,Cosma1,Olivares,Pop}.	
A transportation problem similar to the present  two-phase transportation  problem  was  studied by Sharma et al. \cite{EJOR} with the objective of finding a feasible Phase-I transportation schedule so that the corresponding optimal Phase-II  transportation schedule  is such that the sum of the  overall transportation times of both the phases is minimum. They  proposed an iterative solution  algorithm with promising  computational  efficiency. Although, the authors studied this problem with the objective of minimizing the sum of Phase-I transportation time and the corresponding optimal   Phase-II transportation time, the objective in the problem could also be considered as minimization of the  sum of Phase-I transportation  time and the corresponding Phase-II transportation time instead of corresponding `optimal' Phase-II transportation time. By making this change, the ultimate goal of the problem does not change.   Keeping in view, the  situations where the  budget is as important as the in-time delivery of the product,    the idea of the present study
was conceived (i.e., to discusses the two-phase transportation problem similar to the one discussed in   \cite{EJOR}, as a bi-objective transportation problem  with two objectives viz., sum of transportation times and sum of transportation costs of both the phases). Therefore, it is named as a bi-objective two-phase transportation problem (BTPTP). 
This problem   has its relevance in many real life situations, specifically, in  `fast moving consumable goods' industries where the destinations range from  big cities to small towns and villages. The industry targets on transporting the quantity of its product equal to the market demand in the minimum possible  time as well cost. This is achieved in two phases as Phase-I transportation is done among the preferred links and Phase-II transportation is done among the non-preferred links. The preference to various source-destination links is based on many factors such as higher demand, better routes or lesser taxes etc. After the transportation in Phase-I is complete, the industry would come to know the new outstanding demand of the destinations, and then the non-preferred links are considered for transportation. \\ Another very important situation that gives rise to  the two phase transportation problem   is related to  the natural calamity hit hilly areas  specifically, in which either the rescue teams or the products  like food packets, medicines, blankets etc. are to be transported. The prime objective in such situations  usually,  is to execute the transportation in the   minimum possible time and budget, however, some  other objectives may also be relevant to study. Assume that there are $m$ number of sources from which the $n$ destinations have to be catered. Due to heavy snowfall,  landslide, ongoing construction or for  other such reasons, some of the source destination link roads  may be  closed and therefore, the  air transportation is the only option  left to cater these links. But to arrange the air transportation all of a sudden  in such areas may take quite long.  So,  the administration decides/prefers to   fulfill as much as possible  demand  of the destinations  firstly by transporting along   the open link roads using local  road  transportation.	After the local transport reaches some of the  destinations, the administration  comes to know about the leftover demand of all the destinations  and then  transports the enough amount to these destinations by air  and satisfy their demands. This road and air transportation together,  may be  termed as one  trip  of the  mission in which the demand of the destinations for  a particular time interval is satisfied.	The administration's  objective is to carry out each trip of the mission  in the minimum possible     time as well as cost. Due to the conflicting nature of the two objectives,  it may not be possible to find a single optimal solution, therefore, an iterative algorithm is developed that records all the non-dominated points i.e., the pairs of sum of transportation times and sum of transportation costs of both the phases. No additional subjective preference  is assumed, so all the non-dominated    points are considered equally good.
\par This paper is organized as follows:   Section 2 provides the mathematical formulation of problem BTPTP. An algorithm is developed in  Section 3 along with various theoretical results justifying its working. Numerical illustration is given in Section 4.    Section 5 contains the  computational behavior of the BTPTP-Algorithm for BTPTP instances of different sizes. The concluding remarks are given in Section 6.
\begin{table}[h]\centering                    \caption{Literature Review of Multi-stage and two-stage transportation problems }
	\begin{tabular}{lcc }\label{tt}
		\textbf{Author}&\textbf{Type of  } & \textbf{objectives} \\
		&\textbf{ transportation  problem} &  \textbf{(Minimizing)} \\\hline
		Osman and  & Multi-stage   & Transportation cost \\
		Ellaimony \cite{G} &     &  \\\hline
		Osman  &   Bi-objective, Multi-stage   & Transportation and \\
		Ellaimony \cite{osman,G1} &      &   Deterioration cost\\\hline
		Istranikova \cite{F}&Two-stage   &Transportation cost \\\hline
		Gen et al. \cite{Gen2006} &Two-stage  & Logistics cost  \\\hline
		Ritha et al. \cite{K} & Fuzzy with multi  & Transportation  cost \\
		& -objective constraints &   \\\hline
		Murad et al. \cite{101} & Bi-objective,   & Transportation and  \\
		&  Multistage  &  Deterioration cost \\\hline
		
		Cosma et al. \cite{Cosma} &  Two-stage  &  Fixed costs \\\hline
		Cosma et al. \cite{Cosma1}  & Two-stage   &  Fixed charges and   \\
		&   &   transportation costs   \\\hline
		Calvete et al. \cite{Calvete,Calvete1},   &  Two-stage  &  Fixed costs and  variable \\
		Pop et al. \cite{Pop} &     & costs proportional to the   \\
		&     &   number of units  transported\\ \hline
		Olivares et al.  \cite{Olivares}   & Two-stage, bi-objective& Transportation time, cost  \\\hline
		Sharma et al. \cite{EJOR}& Two-level  &  Sum of transportation \\
		&    &  times of both levels  \\
		&    &     \\
		\hline
		Kaushal and Arora \cite{Kaushal} &Two-stage& Sum of transportation \\
		&    &  times of both levels  \\
		\hline
		Ekta et al. \cite{ekta3}&Three-phase  &  Sum of transportation \\
		&    &  times of all phases  \\
		
		\hline
		Present problem &  Bi-objective, Two-phase &  Sum of transportation times \\
		&    & and  costs of both phases\\\hline
\end{tabular} \end{table}
\section{Mathematical Formulation of BTPTP}	
\begin{center} \begin{figure}[h]
		\begin{minipage}{6.0in}
			\includegraphics[width=\textwidth]{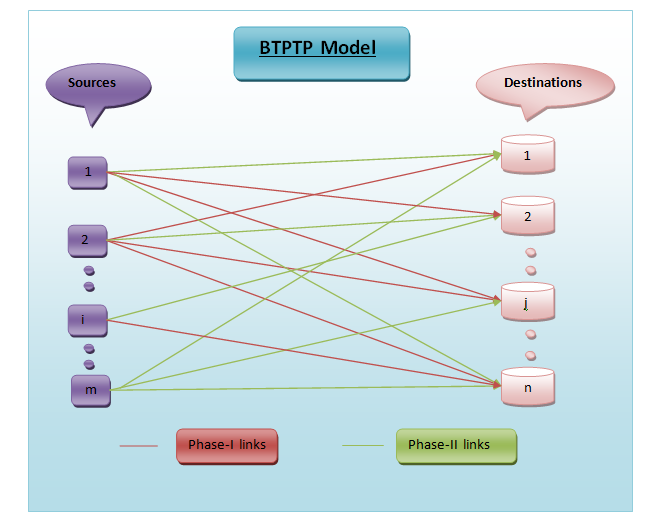}
		\end{minipage}
		{\footnotesize\caption{A  model representing a two phase transportation problem}\label{Model}}
	\end{figure}
\end{center}

Let $I=\{1,2,\dots,m\}$ and $J=\{1,2,\dots,n\}$ represent the set of sources and destinations, respectively, of a balanced  transportation problem with $m$ sources and $n$ destinations. The quantities $a_i$, for $i \in I $ and $b_j$,   for  $j\in J$ denote the availability and demand of the homogeneous product at $i^{th}$ source and $j^{th}$ destination, respectively. The quantity $t_{ij}$ denotes the time consumed in transporting the homogeneous product from $i^{th}$ source to $j^{th}$ destination whereas  $c_{ij}$ denotes the per unit transportation cost.
It is important to note here that 	unlike $c_{ij}$, $t_{ij}$  is independent of the quantity being transported.
\\Further, let $y_{ij}$ and $z_{ij}$ denote the amounts of   product being transported from $i^{th}$ source to $j^{th}$ destination. 
In a balanced two-phase transportation problem, the set $I \times J$ of all the source-destination links  is partitioned into two non-empty disjoint sets $PH_1$ and $PH_2$.    The set $PH_1$ consists of  links  to be used in transportation (of the homogeneous product) during Phase-I and the set $PH_2$  consists of links to be used in transportation during Phase-II.  
In Phase-I of the problem, transportation is carried out ``in parallel" only among the links belonging to $PH_1$    without exhausting  availability and demand of various sources and destinations, respectively. 
After the  transportation in Phase-I of the problem is complete, the leftover available amount of the product at  various sources is transported (in parallel) to satisfy the leftover demands of various destinations by using the links belonging to $PH_2$ only. During each phase of the  transportation, transporting vehicles start simultaneously from all the permissible sources. The objective in each of problems is to minimize the corresponding  transportation time and transportation cost, simultaneously. However, the objective of  the two-phase transportation problem, is to minimize the sum of transportation times and sum of transportation costs, simultaneously. A pictorial layout of this  problem with $m$ sources and $n$ destinations   is given in Figure \ref{Model}. The source-destination links represented by red and green arrows belong to  Phase-I and Phase-II sets, respectively.	 
Now, the bi-objective transportation problems corresponding to Phase-I  and Phase-II   can be  defined  as follows.\\
\textbf{ Phase-I problem }
$$ \min_{Y \in S_{PH_1} }\left\{T_1(Y), C_1(Y)\right\}$$
where $$\begin{array}{lr} T_1(Y)= \displaystyle\max_{(i,j)\in I \times J} \{t_{ij}(y_{ij})\};& C_1(Y)=\displaystyle\sum_{(i,j)} \displaystyle\sum_{\in {I \times J}}c_{ij}y_{ij} \end{array} $$ and
$$    S_{PH_1}=\left\{\begin{array}{cc}
	Y=\{y_{ij}\}_{I \times J} \in R^{m \times n}  &   \left | \begin{array}{l}
		\displaystyle\sum_{j \in J}y_{ij} \leq a_{i},  ~\forall ~ i \in I\\
		\displaystyle\sum_{i \in I} y_{ij} \leq b_{j}, ~\forall ~j \in J\\
		y_{ij} \geq 0, ~\forall~ (i,j) \in PH_1,~y_{ij}=0,~\forall~ (i,j) \in PH_2\\
		a_i,b_j,y_{ij} \mbox{~are~integers and} \displaystyle\sum_i\displaystyle\sum_jy_{ij} \neq 0\\
	\end{array} \right . \end{array} \right\}$$
is the solution set. Here,
$      t_{ij}(y_{ij})=\left\{\begin{array}{ll}
	t_{ij}               &  \mbox{if} ~~y_{ij} >0\\
	0 &  \mbox{otherwise}\end{array}\right.. $\vspace{1mm}\\
\textbf{Phase-II problem corresponding to a feasible solution $Y$ of the Phase-I problem.}
$$ \displaystyle\min_{Z \in S_{PH_2}(Y) }\left\{T_2(Z), C_2(Z)\right\}$$
where $$\begin{array}{lr}T_2(Z)= \displaystyle\max_{(i,j)\in I \times J} \{t_{ij}(z_{ij})\}; &C_2(Z)=\displaystyle\sum_{(i,j) } \displaystyle\sum_{\in I \times J}c_{ij}z_{ij} \end{array} $$
and
$$    S_{PH_2}(Y)=\left\{\begin{array}{cc}
	Z=\{z_{ij}\}_{I \times J} \in R^{m \times n}  &   \left | \begin{array}{l}
		\displaystyle\sum_{j\in J}z_{ij} =  a_{i}-a_i^{'},  ~\forall ~ i \in I\\
		
		\displaystyle\sum_{i \in I} z_{ij} = b_{j}-b_{j}^{'}, ~\forall ~j \in J\\
		
		z_{ij} \geq 0, ~\forall~ (i,j) \in PH_2,~z_{ij}=0,~\forall~ (i,j) \in PH_1\\
		
		z_{ij} \mbox{~are~integers~and} \displaystyle\sum_i\displaystyle\sum_jz_{ij} \neq 0
		
	\end{array} \right . \end{array} \right\}$$
is the solution set. Here $a_{i}^{'}=\displaystyle\sum_{j \in J} y_{ij}$ is the quantity transported from $i^{th}$ source, $b_{j}^{'}=\displaystyle\sum_{i \in I} y_{ij}$ is the quantity transported to $j^{th}$ destination in   Phase-I of the problem and
$      t_{ij}(z_{ij})=\left\{\begin{array}{ll}
	t_{ij}               &  \mbox{if} ~~z_{ij} >0\\
	0 &  \mbox{otherwise}\end{array}\right.
$.  \vspace{3mm}\\
\noindent Now, the bi-objective two-phase transportation problem (BTPTP) can mathematically,  be stated as
$$~~~~~~~~~~~~~~~~~~~~~~~~~~~~~~~~~\min_{Y \in S_{PH_1}, Z \in S_{PH_2}(Y) }\left\{T_{1}(Y)+ T_{2}(Z), C_{1}(Y)+C_{2}(Z)\right\}~~~~~~~~~~~~~~~~~(BTPTP)
$$
In the formulation of the bi-objective transportation problems corresponding to both the phases and hence, in that of the problem BTPTP,  $c_{ij}$ and $t_{ij}$ are assumed to be non-negative, for all $ (i,j)\in I \times J$.\\  \\
A feasible solution $(Y,Z)\in S_{PH_1} \times S_{PH_2}(Y) $ of the problem BTPTP can  be obtained from a feasible solution $X$, of  an equivalent balanced transportation problem $P_0$  defined as follows:
$$ ~~~~~~~~~~~~~~~~~~~~~~~~~~~~~~~~~~~~~~~~~\min_{X \in S }\left\{S_T(X), C(X)\right\}~~~~~~~~~~~~~~~~~~~~~~~~~~~~~~~~~~~~~~~~(P_0)$$   where $S$ is the solution set $S$   defined as
$$
S=\left\{\begin{array}{cc}
	X=\{x_{ij}\}_{I \times J}  \in R^{m \times n}  &   \left | \begin{array}{l}
		\displaystyle
		\sum_{j \in J}x_{ij}=a_{i},  ~\forall ~ i \in I\\
		
		\displaystyle\sum_{i \in I} x_{ij}=b_{j}, ~\forall ~j \in J\\
		x_{ij} \geq 0 \mbox{~and integers},~
		\forall~ i \in I, ~j \in J
	\end{array} \right . \end{array} \right\}                       $$
and $C(X)$ and $S_T(X)$ denote the total transportation cost and the sum of transportation times of Phase-I and Phase-II, respectively, corresponding to an  $X\in S$. Mathematically,
$$\begin{array}{lr}C(X)=\displaystyle\sum_{(i,j)}\displaystyle\sum_{ \in {I \times J}}c_{ij}x_{ij};&~~S_T(X)= \displaystyle\max_{(i,j)\in PH_1} \{t_{ij}(x_{ij})\}
	+\displaystyle\max_{(i,j)\in PH_2} \{t_{ij}(x_{ij})\}
\end{array}$$ and
$     t_{ij}(x_{ij})=\left\{\begin{array}{ll}
	t_{ij}               &  \mbox{if} ~~x_{ij} >0\\
	0 &  \mbox{otherwise}\end{array}\right.
.$\\
Here,  $C(X)$ is a linear function whereas $S_T(X)$ is a concave function \cite{bansal}.\\ \\
Corresponding to a feasible solution $X \in S$ of the problem $P_0$,  a feasible solution $(Y,Z)$ of BTPTP can be defined as   
$$   \begin{array}{lr}
	{ y_{ij}=\left\{ \begin{array}{ll}
			x_{ij}  &    {\rm if}~(i,j)\in PH_1,\\
			0       &  {\rm if}~(i,j)\in PH_2                           \end{array} \right .}~ \mbox{and}&{z_{ij}=\left\{ \begin{array}{ll}
			0   &   {\rm if} ~ (i,j)\in PH_1,\\
			x_{ij}  & {\rm if}~(i,j)\in PH_2                           \end{array} \right .}
\end{array} $$
Equivalently,  \begin{equation}\label{15}X=Y +Z\end{equation}
Clearly, $Y \in S_{PH_1}$ and $Z \in S_{PH_2}(Y)$.\\
Moreover,
$$  \begin{array}{ll}
	S_T(X)  &=\displaystyle\max_{ (i,j)\in PH_1} \{t_{ij}(x_{ij})\}
	+\displaystyle\max_{ (i,j)\in PH_2} \{t_{ij}(x_{ij})\}\\
	&= \displaystyle\max_{(i,j)\in I \times J } \{t_{ij}(y_{ij})\}     +\displaystyle\max_{(i,j)\in I \times J} \{t_{ij}(z_{ij})\}\\&=T_1(Y)+T_2(Z)
\end{array}$$ and
$$  \begin{array}{ll}
	C(X)  &=\displaystyle\sum_{(i,j) }\displaystyle\sum_{\in {I \times J}}c_{ij}x_{ij}\\&=
	\displaystyle\sum_{(i,j)}\displaystyle\sum_{ \in {PH_1}}c_{ij}x_{ij}+\displaystyle\sum_{(i,j)}\displaystyle\sum_{ \in {PH_2}}c_{ij}x_{ij}\\
	&=\displaystyle\sum_{(i,j)}\displaystyle\sum_{ \in {I \times J}}c_{ij}y_{ij}+ \displaystyle\sum_{(i,j)}\displaystyle\sum_{\in {I \times J}}c_{ij}z_{ij}\\&=C_1(Y)+C_2(Z)
\end{array}$$
It shows that a feasible solution $(Y,Z) $ of BTPTP can be obtained from 
a feasible solution $X$ of the problem $P_0$. \\
Therefore, the proposed technique solves the  problem  $P_0$ and its time restricted versions to find all the non-dominated pairs of the problem BTPTP. 
\subsection{Definitions}
In this sub-section, we provide various general  definitions \cite{Kaur} in the context of   present problem as follows.   
\begin{definition} \textbf{{Feasible pair}}. A pair $(S_T, C) $ of sum of transportation times and total  transportation cost corresponding to the problem $P_0$  is called a feasible pair if $\exists$ a feasible solution $X \in S$ such that  $(S_T, C)=(S_T(X), C(X))$.
\end{definition}
\begin{definition}
	\textbf{{Pareto optimal  solution}}. A solution $\hat{X} \in S$  is called a Pareto optimal solution of the problem $P_0$ if there does not exist a feasible solution $X \in S $ such that for the pairs  $(S_T(X), C(X))$ and $(S_T(\hat{X}), C(\hat{X}))$, either $S_T(X)<S_T(\hat{X})$ and $C(X)\leq C(\hat{X})$  or $S_T(X) \leq S_T(\hat{X})$ and $C(X)<  C(\hat{X})$, holds. 	 \vspace{1mm}\\ 
	Further, a pair $(S_T(\hat{X}), C(\hat{X}))$ corresponding to the Pareto optimal solution $\hat{X}$  will be called a non-dominated  point of the problem $P_0$.\end{definition}
\begin{definition} \textbf{Dominated solution}. A solution $\bar{X}\in S$  is called a dominated solution of the problem $P_0$ if $\exists$ a feasible solution $X \in S $ such that for the pairs  $(S_T(X), C(X))$ and $(S_T(\bar{X}), C(\bar{X}))$, either $S_T(X)<S_T(\bar{X})$ and $C(X) \leq C(\bar{X})$  or $S_T(X)\leq S_T(\bar{X})$ and $C(X)<  C(\bar{X})$, holds. \end{definition}
\section{Theoretical Development}
\noindent To develop a solution technique for the problem BTPTP that records all of its non-dominated  solutions, we initially formulate and solve  a cost minimizing transportation problem (CMTP) named as  $CP_0$, associated to the problem $P_0$. Mathematically, the problem $CP_0$ is defined as  $$
~~~~~~~~~~~~~~~~~~~~~~~~~~~~~~\min_{X\in S}C(X)=\min_{X\in S} \displaystyle\sum_{(i,j)} \displaystyle\sum_{\in {I \times J}}c_{ij}x_{ij}~~~~~~~~~~~~~~~~~~~~~~~~~~~~~~~~~~~~~~ (CP_0)$$
In this problem, the objective  is to minimize the total transportation cost associated to the problem $P_0$ and the solution set is same as that of the problem $P_0$.\\
From an OFS $X^0=\{x^0_{ij}\}_{I \times J}$ say, of the problem $CP_0$, the optimal  transportation cost $C_0$ is recorded as $\displaystyle\sum_{(i,j)} \displaystyle\sum_{\in {I \times J}}c_{ij}x^0_{ij}$.\\
Using relation $(1)$, define the corresponding feasible solutions of Phase-I and Phase-II problem  as $Y^0=\{y^0_{ij}\}_{I \times J} $ and $Z^0 =\{z^0_{ij}\}_{I \times J}$  yielding  the  transportation  times as  $T^0_1$  and $T^0_2$, respectively, as  
$$\begin{array}{ccc}
	& T^0_1= \displaystyle\max_{(i,j)\in I \times J} \{t_{ij}(y^0_{ij})\}~~~~~\mbox{and}~~~~&T^0_2= \displaystyle\max_{(i,j)\in I \times J} \{t_{ij}(z^0_{ij})\}\\
\end{array} $$
where the transportation times $t_{ij},~\forall~(i,j)\in I \times J$ are same as that of the problem $P_0$. Also, $$C_0=\displaystyle\sum_{(i,j)} \displaystyle\sum_{\in {I \times J}}c_{ij}x^0_{ij}=\displaystyle\sum_{(i,j)} \displaystyle\sum_{\in {I \times J}}c_{ij}y^0_{ij}+\displaystyle\sum_{(i,j)} \displaystyle\sum_{\in {I \times J}}c_{ij}z^0_{ij}$$
Let $S_{T_0}=T^0_1+T^0_2$. 	 The pair $(S_{T_0}, C_0)$ is considered as the first feasible pair of the problem BTPTP which consists of the optimal total transportation cost and  sum of the corresponding  transportation times of Phase-I and Phase-II.
Now, various time restricted versions of the problem $P_0$ are solved to check the Pareto optimality of the pair $(S_{T_0}, C_0)$.
We call these  restricted CMTPs corresponding to a particular sum of transportation times as  ``sub-problems restricted at the sum $S_{T_0}$". Mathematically, the sub-problems  restricted at $S_{T_0}$ are defined as follows:  
\begin{itemize}
	\item{ Sub-problem $P_{[T^0_1,T^0_2)}^1(S_{T_0})$
		
		$$\begin{array}{lr}
			\displaystyle\min_{X \in S}\displaystyle\sum_{i \in I}\displaystyle\sum_{j \in J}c'_{ij}x_{ij}	&~\mbox{where}~{c'_{ij}:=\left\{ \begin{array}{lcr}
					c_{ij}  &   {\rm if~ for } ~ (i,j)\in PH_1,& t_{ij} \leq T^0_1\\
					M                  &   {\rm if~ for }~(i,j)\in PH_1,&                          t_{ij} > T^0_1 \\
					c_{ij}  &   {\rm if~ for } ~ (i,j)\in PH_2,& t_{ij} < T^0_2\\
					M                  &   {\rm if~ for }~(i,j)\in PH_2,&                          t_{ij} \geq T^0_2
				\end{array} \right. }
		\end{array}$$			
		where $M$ is a very large positive number.}

	\item{Sub-problem $P_{(T^0_1,T^0_2]}^2(S_{T_0})$
		$$\begin{array}{lr}
			\displaystyle\min_{X \in S}\displaystyle\sum_{i \in I}\displaystyle\sum_{j \in J}c'_{ij}x_{ij}	&~\mbox{where}~{c'_{ij}:=\left\{ \begin{array}{lcr}
					c_{ij}  &   {\rm if~ for } ~ (i,j)\in PH_1,& t_{ij} < T^0_1\\
					M                  &   {\rm if~ for }~(i,j)\in PH_1,&                          t_{ij} \geq T^0_1 \\
					c_{ij}  &   {\rm if~ for } ~ (i,j)\in PH_2,& t_{ij} \leq T^0_2\\
					M              &   {\rm if~ for }~(i,j)\in PH_2,&                          t_{ij} > T^0_2
				\end{array} \right.}
		\end{array}$$}
\end{itemize}
Let $(T_{1})_{{max}}=t_1^\alpha$ be the largest time entry in the time matrix corresponding to the set $PH_1$. Calculate $T'_1=\min\{S_{T_0}, (T_{1})_{{max}}\}$.
Further, arrange  all the time entries  of the time matrix in the increasing order corresponding to the links in $PH_1$ that lies between $T^0_1$ and $T'_1$ as follows:	$$T^0_1 < T_{1.1}< T_{1.2}<\cdots< T_{1.r}<T'_1$$.
Suppose, there are $r$ such time entries  then, the corresponding  $(r+1)$ sub-problems will be  defined as 
\begin{itemize}
	\item{Sub-problems $P_{[T_{1.p},S_{T_0}- T_{1.p})}^{3.p}(S_{T_0}),~p=1,2,\dots,r+1.$

		$$\begin{array}{lr}
			
			\displaystyle\min_{X \in S}\sum_{i \in I}\sum_{j \in J}c'_{ij}x_{ij}&~\mbox{where}~
			c'_{ij}:=\left\{ \begin{array}{lcc}
				c_{ij}  &   {\rm if~ for } ~ (i,j)\in PH_1,& t_{ij} \leq T_{1.p}\\
				M                  &   {\rm if~ for }~(i,j)\in PH_1,&                          t_{ij} >  T_{1.p} \\
				c_{ij}  &   {\rm if~ for } ~ (i,j)\in PH_2,& t_{ij} < S_{T_0}- T_{1.p}\\
				M                  &   {\rm if~ for }~(i,j)\in PH_2,&                          t_{ij} \geq S_{T_0}- T_{1.p}
			\end{array} \right. \end{array}$$ Clearly, $T_{1.r+1}= T'_1 $\vspace{2mm}\\      
		\begin{remark}
			For $p=r+1$, the sub-problem $P_{[T_{1.p},S_{T_0}- T_{1.p})}^{3.p}(S_{T_0})$ becomes
			$$\begin{array}{lr}
				\displaystyle\min_{X \in S}\sum_{i \in I}\sum_{j \in J}c'_{ij}x_{ij}&~\mbox{where}~
				c'_{ij}:=\left\{ \begin{array}{lcc}
					c_{ij}  &   {\rm if~ for } ~ (i,j)\in PH_1,& t_{ij} \leq T'_1\\
					M                  &   {\rm if~ for }(i,j)\in PH_1,&                          t_{ij} >  T'_1\\
					c_{ij}  &   {\rm if~ for } ~ (i,j)\in PH_2,& t_{ij} < S_{T_0}-T'_1\\
					M                  &   {\rm if~ for }(i,j)\in PH_2,&                          t_{ij} \geq S_{T_0}-T'_1
				\end{array} \right. \end{array} $$
			Depending upon the value of $T'_1$, there are either $r$ problems (for $T'_1=S_{T_0}$) or $r+1$ sub-problems (for $T'_1= (T_1)_{max}$). For first case, due to the blocking of all links in Phase-II, no transportation can take place. Therefore, a feasible solution of this solution can not be obtained. 			
\end{remark}}     	\end{itemize}

Let $(T_{2})_{{max}}=t_1^\beta$ be the largest time entry in the time matrix corresponding to the set $PH_2$. Calculate $T'_2=\min\{S_{T_0}, (T_{2})_{{max}}\}$.
Further, arrange  all the time entries  of the time matrix in the increasing order corresponding to the links in $PH_2$ that lie between $T^0_2$ and $T'_2$ as follows:	$$T^0_2 < T_{2.1}< T_{2.2}<\cdots< T_{2.s}<T'_2$$.
Suppose, there are $s$ such time entries  then, the corresponding  $(s+1)$ sub-problems will be  defined as 
\begin{itemize}
	\item{Sub-problems $P_{(S_{T_0}-T_{2.q},T_{2.q}]}^{4.q}(S_{T_0}),~q=1,2,\dots,s+1.$
		$$\begin{array}{lr}
			\min_{X \in S}\sum_{i \in I}\sum_{j \in J}c'_{ij}x_{ij}&~\mbox{where}~
			c'_{ij}:=\left\{ \begin{array}{lcc}
				c_{ij}  &   {\rm if~ for } ~ (i,j)\in PH_2,& t_{ij} \leq T_{2.q}\\
				M  &   {\rm if~ for } ~ (i,j)\in PH_2,& t_{ij} > T_{2.q}\\
				c_{ij}  &   {\rm if~ for } ~ (i,j)\in PH_1,& t_{ij} <S_{T_0}- T_{2.q}\\
				M                  &   {\rm if~ for }~(i,j)\in PH_1,&                          t_{ij} \geq  S_{T_0}-T_{2.q} \\
			\end{array} \right. \end{array}$$ 
		Clearly, $T_{2.s+1}= T'_2 $.\\
		Remark 1 for the sub-problem
		$P_{[T_{1.p},S_{T_0}- T_{1.p})}^{3.p}(S_{T_0})$ for $p=r+1$, also holds for the problem $P_{(S_{T_0}- T_{2.q},T_{2.q}]}^{4.q}(S_{T_0})$ for $q=s+1$.}
\end{itemize}
\begin{remark}
	Clearly, the sub-problems $P_{[T^0_1,T^0_2)}^1(S_{T_0})$ and  $P_{[T_{1.p},S_{T_0}- T_{1.p})}^{3.p}(S_{T_0}),~p=1,2,...,r+1$ are of the form  $P_{[a^0_i,b^0_i)}^{i}(S_{T_0})$~; for $i=1,3.p$, for some $a^0_i,~b^0_i \in \mathbf{R}^+$ whereas the sub-problems $P_{(T^0_1,T^0_2]}^2(S_{T_0})$ and $P_{(S_{T_0}- T_{2.q},T_{2.q}]}^{4.q}(S_{T_0}),~q=1,2,...,s+1$
	are of the form $P_{(a^0_i,b^0_i]}^{i}(S_{T_0})$~; for $i=2,4.q$. \end{remark} 
\begin{remark}\label{mee} In a sub-problem of the form $P_{[a^0_i,b^0_i)}^{i}(S_{T_0})$,  $c'_{ij}=c_{ij}$, for all $(i,j)\in PH_1$ with corresponding $t_{ij} \leq a^0_i$ and $c'_{ij}=M$, for all $(i,j)\in PH_1$ with corresponding $t_{ij}> a^0_i$. Similarly, $c'_{ij}=c_{ij}$, for all $(i,j)\in PH_2$ with corresponding $t_{ij} < b^0_i$ and $c'_{ij}=M$, for all $(i,j)\in PH_2$ with corresponding $t_{ij}\geq b^0_i$.  Therefore,  the numbers  $a^0_i$ and  $b^0_i$ act as the upper bounds for the time entries of the links of Phase-I and Phase-II sets, respectively.
	Clearly, the parenthesis `[' and `)' indicate about which of the upper bound is
	to be included or excluded, respectively, while blocking the links. \vspace{1mm}\\
	Similar argument hold for the problem 	$P_{(a^0_i,b^0_i]}^{i}(S_{T_0})$ as well. \\
	Since, all the source-destination links in the problem $CP_0$ possess their original transportation costs, therefore, in this problem, the upper bounds for  the time entries of the links of Phase-I and Phase-II sets may be taken  as infinite. Equivalently, the problem $CP_0$ may be written as $P^i_{(\infty, \infty)}$, for any $i\in \{1,2,3.p,4.q\}$.  \end{remark}
\begin{remark} A sub-problem of the form  $P_{[a^0_i,b_i^0)}^{i}(S_{T_0})$ is equivalent to a sub-problem of the form  $P_{[a_i^{0'},b_i^{0'}]}^{i}(S_{T_0})$, for some $a_i^{0'},b_i^{0'}\in \mathbf{R}^+$ with   $a_i^{0'}=a^0_i$ and $b_i^{0'}<b_i^0$ such that there is no time entry between $b_i^{0'}$ and $b_i^0$. This argument is true for the sub-problems of the form $P_{(a_i^0,b_i^0]}^{i}(S_{T_0})$ also, if  $a_i^{0'}<a_i^0$ such that there is no time entry between $a_i^{0'}$ and $a_i^0$ and $b_i^{0'}=b_i^0$.  Therefore, we may, in general,   denote all the above sub-problems as $P_{[a_i^0,b_i^0]}^{i}(S_{T_0})$ also, for  some $a_i^0,b_i^0 \in \bf{R}^+$.
\end{remark}	 
\begin{remark}\label{rr}
	Since, the above sub-problems are  restricted versions of the CMTP $CP_0$ only, therefore, these sub-problems can be solved by using any solution technique available in literature, for solving any CMTPs.
\end{remark}
\begin{remark} \label{remi}
	To understand the formulation of various sub-problems restricted at a particular sum, let us consider a randomly generated BTPTP instance given in Table \ref{tabexa}, named as  $3 \times 5$ BTPTP 1, in which the cells with bold entries denote the Phase-I links while the others denote Phase-II links. Also, at each source-destination link, the number on top left corner denotes the transportation cost $c_{ij},~\mbox{for} ~i \in I,~ j\in J$ whereas the number on bottom right corner denotes the transportation time $t_{ij}$.  The quantities $a_i$ and $b_j$ denote the availability at $i^{th}$ source and demand at $j^{th}$ destination, respectively.\\
	\noindent The distinct time entries corresponding to  Phase-I and Phase-II links of the given BTPTP are given below: \vspace{1mm}\\Time entries of Phase-I links:~~~$2<3<4<5<6<7<9$.\\  ~~~~~~~~~~~~~~~~Time entries of Phase-II links:~ $3<4<6<8<9<12$.\vspace{2mm}\\
	\begin{table}[h]
		\caption{A $3 \times 5$  BTPTP 1}
		\centering
		\begin{tabular}{c||c|c|c|c|c||c}  \label{tabexa}            & $D_1$ & $D_2$ & $D_3$ & $D_4$ & $D_5$ & $a_i~\downarrow$   \\
			\hline $S_1$ & \begin{tabular}{lr}
				$ 10$&\\& $3$
			\end{tabular}   & \begin{tabular}{lr}
				{$\mathbf{2}$} &\\ &\bf\textbf{3}
			\end{tabular}   &\begin{tabular}{lr}
				{$\mathbf{5}$} &\\ &\bf\textbf{7}
			\end{tabular}   &\begin{tabular}{lr}
				$9$ &\\ &$4$
			\end{tabular}   & \begin{tabular}{lr}
				$8$ &\\ &$6$
			\end{tabular}   &  7\\
			\hline
			$ S_2$   & \begin{tabular}{lr}
				{$\mathbf{8}$} &\\ &\bf\textbf{2}
			\end{tabular}   &  \begin{tabular}{lr}
				$4$ &\\ &$6$
			\end{tabular}     &  \begin{tabular}{lr}
				{$\mathbf{11}$} &\\ &\bf\textbf{4}
			\end{tabular}   & \begin{tabular}{lr}
				$7$ &\\ &$9$
			\end{tabular}   &  \begin{tabular}{lr}
				{$\mathbf{6}$} &\\ &\bf\textbf{9}
			\end{tabular}    & 8   \\
			\hline
			$ S_3$ & \begin{tabular}{lr}
				$3$ &\\&$12$
			\end{tabular}     & \begin{tabular}{lr}
				{$\mathbf{5}$} &\\ &\bf\textbf{6}
			\end{tabular}    & \begin{tabular}{lr}
				$4$ &\\ &$8$
			\end{tabular}   & \begin{tabular}{lr}
				{$\mathbf{6}$} &\\ &\bf\textbf{5}
			\end{tabular}   & \begin{tabular}{lr}
				$9$ &\\ &$9$
			\end{tabular}    &9\\
			\hline $b_j~\longrightarrow$ & 3 & 5 & 6 & 6 & 4 &
		\end{tabular}
	\end{table}
	\noindent One of the feasible solutions of this problem is  obtained by using north-west corner rule for solving  its corresponding CMTP (i.e., the problem $CP_0$ in which only the costs of all the links  are considered). This OFS is given in  Table \ref{tabexa2}.  The encircled numbers  in some of the cells  in   Table \ref{tabexa2}, denote the quantity of the homogeneous product to be transported along that link.  Further,  if we consider  this feasible solution to be $X=\{x_{ij}\}_{I \times J}\in S$, the corresponding Phase-I and Phase-II transportation times are noted as $\displaystyle\max_{(i,j)\in PH_1} \{t_{ij}(x_{ij})\}=5$ and $
	\displaystyle\max_{(i,j)\in PH_2} \{t_{ij}(x_{ij})\}=9$ while the corresponding sum of transportation costs of both the phases is noted as $\displaystyle\sum_{(i,j)}\displaystyle\sum_{ \in {PH_1}}c_{ij}x_{ij}+\displaystyle\sum_{(i,j)}\displaystyle\sum_{ \in {PH_2}}c_{ij}x_{ij}=181$, thereby, yielding a feasible pair of sum of transportation times and sum of transportation costs as $(14, 181)$.
	\begin{table}[h]			
		\caption{A feasible solution of the $3 \times 5$  BTPTP 1}
		\begin{tabular}{c||c|c|c|c|c||c}  \label{tabexa2}            & $D_1$ & $D_2$ & $D_3$ & $D_4$ & $D_5$ & $a_i$   \\
			\hline $S_1$ & \begin{tabular}{lll}
				$10$&&\\&\textcircled{3}&\\&& $3$
			\end{tabular}   & \begin{tabular}{lll}
				{$\mathbf{2}$} &&\\ &\textcircled{4}&\\&&\bf\textbf{3}
			\end{tabular}   &\begin{tabular}{lll}
				{$\mathbf{5}$} &&\\ &&\\&&\bf\textbf{7}
			\end{tabular}   &\begin{tabular}{lll}
				$9$ &&\\ &&\\&&$4$
			\end{tabular}   & \begin{tabular}{lll}
				$8$ &&\\ &&\\&&$6$
			\end{tabular}   &  7\\
			\hline
			$ S_2$   & \begin{tabular}{lll}
				{$\mathbf{8}$} &&\\ &&\\&&\bf\textbf{2}
			\end{tabular}   &  \begin{tabular}{lll}
				$4$ &&\\ &\textcircled{1}&\\&&$6$
			\end{tabular}     &  \begin{tabular}{lll}
				{$\mathbf{11}$} &&\\ &\textcircled{6}&\\&&\bf\textbf{4}
			\end{tabular}   & \begin{tabular}{lll}
				$7$ &&\\ &\textcircled{1}&\\&&$9$
			\end{tabular}   &  \begin{tabular}{lll}
				{$\mathbf{6}$} &&\\ &&\\&&\bf\textbf{9}
			\end{tabular}    & 8   \\
			\hline
			$ S_3$ & \begin{tabular}{lll}
				$3$ &&\\ &&\\&&$12$
			\end{tabular}     & \begin{tabular}{lll}
				{$\mathbf{5}$} &&\\ &&\\&&\bf\textbf{6}
			\end{tabular}    & \begin{tabular}{lll}
				$4$ &&\\ &&\\&&$8$
			\end{tabular}   & \begin{tabular}{lll}
				{$\mathbf{6}$} &&\\ &\textcircled{5}&\\&&\bf\textbf{5}
			\end{tabular}   & \begin{tabular}{lll}
				$9$ &&\\ &\textcircled{4}&\\&&$9$
			\end{tabular}    &9\\
			\hline $b_j$ & 3 & 5 & 6 & 6 & 4 &
		\end{tabular}
	\end{table}\\
	Now, various sub-problems restricted at the sum $14$ can be obtained in accordance with those obtained above  corresponding to the sum $S_{T_0}$, as  $P_{[5,9)}^1(14)$, $P_{(5,9]}^2(14)$, $P_{[6,6]}^{3.1}(14)$, $P_{[7,6]}^{3.2}(14)$, $P_{[9,4]}^{3.3}(14)$ and $P_{[2,12]}^{4.1}(14)$.
	as per all possible combinations of the distinct time entries corresponding to the set $PH_1$ and $PH_2$, respectively. 
	These sub-problems are depicted graphically in $(T_1,T_2)$ plane, in Figure \ref{S14i}.	 
	\begin{landscape}
		\begin{figure}\label{S14i}
			\begin{minipage}{6.5in}
				\centering					
				\includegraphics[width=\textwidth]{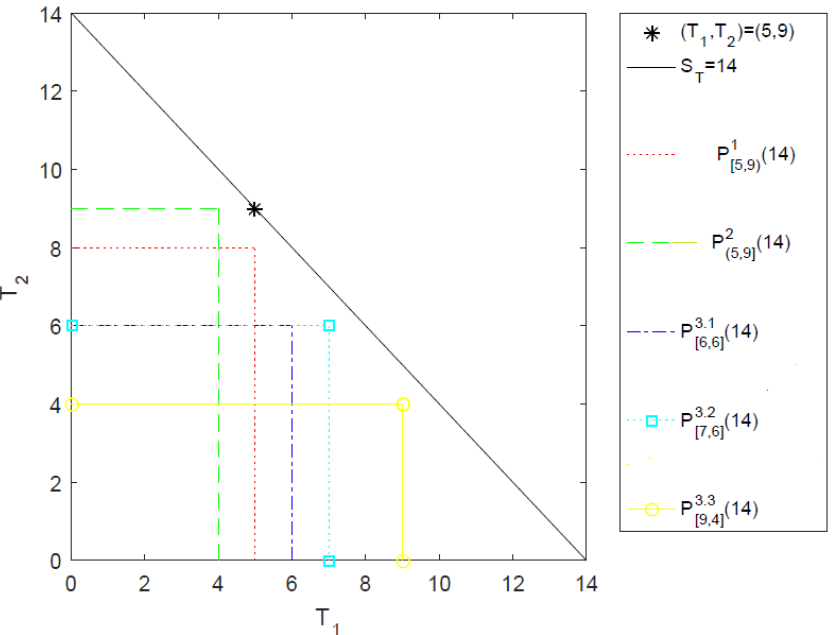}					
			\end{minipage}
			{\footnotesize\caption{Upper bounds for time entries of various sub-problems restricted at the sum 14.}\label{1}}
		\end{figure}
	\end{landscape}
	\begin{remark}\label{ballo}
		It is clear that the sub-problems $P_{[a_i^0,b_i^0]}^{i}(S_{T_0})$, for  all values of $i$, are restricted versions of the CMTP $CP_0$ in which  some links of Phase-I and Phase-II sets  are blocked for transportation by assigning them a very high transportation cost $M$ . These changes in their original costs are made in such a way that a feasible solution of each of these sub-problems (corresponding to which there is no allocation on the blocked links), if exists, would provide a sum of transportation times  strictly less than the sum $S_{T_0}$. We call such a feasible solution  as its M-feasible solution which mathematically, is defined as follows.
	\end{remark}
	\begin{definition} \textbf{ M-feasible solution} \cite{Kaur}. A feasible solution $X=\{x_{ij}\}_{I \times J}\in S$ of a sub-problem, is called an {\bf M-feasible solution}  if $x_{ij}=0,~\forall~(i,j) \in I \times J$ for which $c'_{ij}=M$.\\
		\noindent Further, a sub-problem is called M-feasible if it possesses an M-feasible solution otherwise it is called non M-feasible problem.\\
		An optimal M-feasible solution (OMFS) of a sub-problem is a feasible solution which is optimal as well as M-feasible.\end{definition}
\end{remark}
Thus, it clearly follows  that an OMFS of each of these problems (if exists), will provide a feasible pair of sum of transportation times and sum of transportation costs in which the sum of transportation times is strictly less than $S_{T_0}$.
\begin{remark}
	In the following text, an OMFS of  the sub-problem $P_{[a_i^0,b_i^0]}^{i}(S_{T_{0}})$, for each $i \in \{1,2,3.p,4.q\}$, is denoted by $X^0_{i}$ and the corresponding optimal total transportation  cost    is denoted by \textbf{$CP_{[a_i^0,b_i^0]}^i(S_{T_{0}})$}, i.e.,
	$
	CP_{[a_i^0,b_i^0]}^i(S_{T_{0}})= C(X^0_{i})$,  for all values of $i$.
\end{remark}
\begin{remark}
	As mentioned in Remark \ref{mee}, the sub-problems $P_{[a^0_i,b_i^0]}^{i}(S_{T_0})$, for  all values of $i$, are   restricted versions of the problem $CP_0$, therefore, an M-feasible solution of each of these sub-problems will be
	a feasible solution of the problem $CP_0$. This argument is true for their OMFSs as well.
	Therefore, it follows that $CP_{[a_i^0,b_i^0]}^{i}(S_{T_0})\geq C_0,~\forall ~i$.
\end{remark} 
\begin{definition} (Minor/Major sub-problem).  Given two sub-problems $P_{[a_i,b_i]}^{i}(S_{T})$  and  $P_{[a_j,b_j]}^{j}(S_{T})$ restricted at a sum $S_{T}$, the problem $P_{[a_i,b_i]}^{i}(S_{T})$ is said to be a minor sub-problems of $P_{[a_j,b_j]}^{j}(S_{T})$ if $a_i \leq a_j$ and $b_i \leq b_j$ is satisfied with one of these in-equations being satisfied with strict inequality \cite{Kaur}.

\end{definition}
\begin{remark}\label{2222}
	The major and minor sub-problems have same  feasible region  but different cost matrices for set of open  links  of the latter is contained in that of the former. Therefore, every MFS of a minor sub-problem is also an MFS of major sub-problem but the converse is not true. Also, the non M- feasibility of a major sub-problem leads to the non M-feasibility of all of its minor sub-problems. \end{remark}
The idea of the proposed algorithm is to systematically record non-dominated solutions/points such that the sum of total transportation time decreases  strictly.  For this, the proposed algorithm, at each iteration,  solves various major sub-problems and scans their OMFSs only. Depending upon the optimal cost  of all these major sub-problems  dominance and non-dominance of the in-hand pair is checked. \\
As mentioned above,  out of all the sub-problems $P_{[a_i^0,b_i^0]}^{i}(S_{T_0})$, only the major sub-problems are solved. This is  for the reason explained as follows:\\
Minor and major sub-problems have same feasible region but there is difference in their transportation costs along various  links depending upon the upper bounds on the time entries corresponding to $PH_1$ and $PH_2$ links. Remark 10 implies that every MFS of a minor sub-problem is also an MFS of the major sub-problem, therefore, the set of MFSs of a minor sub-problem is always contained in that of its major sub-problems which implies that the optimal cost of the minor sub-problem would always be greater than or equal to that of the major sub-problem. Since, the proposed algorithm systematically decreases sum of  transportation times whereas the sum of transportation costs  increases, therefore, the minor sub-problem would appear as a major sub-problem in the next iterations, if it has to provide a non-dominated pair.  \\	 
Now,  corresponding to  an OMFS of each of the major sub-problems (if exists), obtain the sum of  transportation times and  the sum of  transportation costs of both the phases (as discussed in Remark \ref{remi}) and yield various feasible pairs.	If the sum of transportation costs corresponding to all of  these feasible pairs are strictly greater than $C_0$, then note the pair $(S_{T_0}, C_0)$ as the first overall non-dominated   pair (see Theorem \ref{wahet1}). Now, construct  a set  consisting of all the other non-dominated  feasible pairs  till now, (obtained from the OMFSs of the major sub-problems restricted at the sum $S_{T_0}$).
Out of these pairs, select a pair with
the minimum   sum of transportation costs. Call this pair as $(S_{T_1}, C_1)$ and construct the sub-problems restricted at $S_{T_1}$ and repeat the above process by replacing $S_{T_0}$ by $S_{T_1}$ in order to find second overall non-dominated pair $(S_{T_2}, C_2)$. Clearly, $S_{T_1}<S_{T_0}$, however, $C_1>C_0$. This process is continued until a pair with minimum sum of transportation times and the corresponding sum of transportation costs is obtained (see Theorem \ref{wahet2}). \\
On the other hand, if the sum of transportation costs corresponding to at least one of the feasible pairs so obtained, is equal to $C_0$, then the pair  $(S_{T_0}, C_0)$ would be dominated by that pair.  This is due to the reason that in this particular pair $(S'_T, C')$ say,  the sum of transportation times  $S'_T$ is strictly less than $S_{T_0}$ but, the sum of transportation costs $C'$ would be equal to $C_0$ (see Theorem \ref{wahel}). Therefore, we replace the pair  $(S_{T_0}, C_0)$ by the pair  $(S'_T, C')$  and
construct
the sub-problems restricted at the sum $S'_T$ . This process is repeated until we obtain the pair  $(S_{T_0}, C_0)$  as the first (overall) non-dominated  pair of sum of transportation times and sum of transportation costs.
\begin{remark}It is worth-mentioning that all the definitions and remarks that hold for the sub-problems restricted at the sum $S_{T_0}$, will also hold for the sub-problems restricted at the  sum $S_{T_1}$, $S_{T_2}$ and so on.\\
	In the following section, the process discussed above is written in the form of an algorithm  named as `BTPTP-Algorithm'.\end{remark}
\subsection{The BTPTP-Algorithm}\label{algorithm}
\begin{landscape}
	\begin{figure}[h]
		\begin{minipage}{8.0in}
			\centering
			\includegraphics[width=\textwidth]{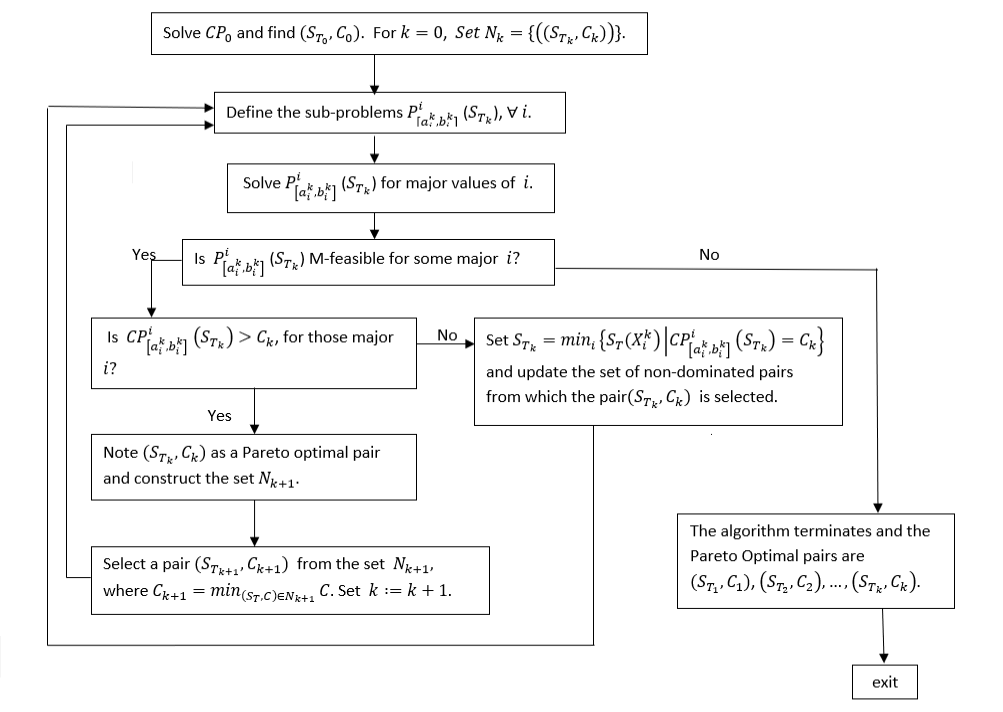}
		\end{minipage}    {\footnotesize\caption{Outline of the BTHTP-algorithm}\label{ttt}}    \end{figure}    \end{landscape}

The working of the algorithm is explained in the context of a feasible pair $(S_T, C)$ of the problem BTPTP.
Before giving the step-wise description of the  BTPTP-Algorithm, following points should be noted. \\
(i)	In the following text, by the term `major value of $i$', we mean `that value of $i$ for which $i^{th}$ sub-problem $P_{[a_i,b_i]}^{i}(S_{T})$ is major'.\\(ii)
If a result holds for all major values of $i$, then for $i=3.p$ (or $4.q$), it means that the result holds for all major values of $p$ (or $q$).\\(iii)
If a result holds at least for one major $i$, then for $i=3.p$ (or $4.q$), it means that the result holds at least for one major values of $p$ (or $q$).
The algorithmic representation of the proposed algorithm is given as follows:\\
\textbf{Step 1.}
\begin{enumerate}
	\item Enter the number of sources, $m$.
	\item Enter the number of destinations, $n$.
	\item Generate the set of transportation links $I \times J$. 
	\item Enter the availability of   sources.
	\item Enter demand of destinations. 
	\item Construct the Phase-I and Phase-II sets namely, $PH_1$ and $PH_2$.	
\end{enumerate}
\textbf{Step 2.}
\begin{enumerate}
	\item Enter the time matrix  $T=\{t_{ij}\}_{m \times n}$.
	\item Enter the time matrix  $C=\{c_{ij}\}_{m \times n}$. 	 
\end{enumerate}
\textbf{Step 3.}
Construct the CMTP $CP_0$. Apply any available solution technique for solving a CMTP and obtain its OFS as $X^{0}$.  From   $X^{0}$, find the optimal cost $C_{0}$ and the corresponding Phase-I and Phase-II transportation times $T^0_{1}$, $T^0_{2}$ and hence, their sum $S_{T_0}$.  Set $k=0$ and define the set $N_k=\{(S_{T_k}, C_k)\}$, the collection of non-dominated pairs till now.\\
\textbf{Step 4.} Define the sub-problems restricted at the sum $S_{T_k}$ i.e., $P_{[a^k_i,b^k_i]}^i(S_{T_k})$, for all values of $i$ and the corresponding values of  $a^k_i$ and $b^k_i$. Now,  solve these sub-problems only for major values of $i$. Let an OMFS of the sub-problem  $P_{[a^k_i,b^k_i]}^i(S_{T_k})$  (if exists), be denoted by  $X^k_{i},\forall~i$.
\begin{enumerate}
	\item {\begin{enumerate}\item If $P_{[a^k_i,b^k_i]}^i(S_{T_k})$ is non M-feasible for all major values of $i$, then  go to the terminal step otherwise go to  (2).\end{enumerate} }
	\item {\begin{enumerate}\item If  $CP_{[a^k_i,b^k_i]}^{i}(S_{T_k}) > C_k ~\mbox{for all those major values of}~i$ for which the sub-problem $P_{[a^k_i,b^k_i]}^i(S_{T_k})$ is M-feasible, then note $(S_{T_k}, C_k)$ as an overall non-dominated   pair and go to (3), otherwise go to (2)(b).
			\item  If $CP_{[a^k_i,b^k_i]}^{i}(S_{T_k}) = C_k~\mbox{ at least for one major}~i$, then the pair $(S_{T}(X^k_{i}), C(X^k_{i}))$ obtained from an OMFS $X^k_{i}$ of the sub-problem $P_{[a^k_i,b^k_i]}^i(S_{T_k})$ for such major values of $i$, would dominate the pair $(S_{T_k}, C_k)$ as $S_{T}(X^k_{i}) < S_{T_k}$, for all such values of $i$. Therefore, update the value of $S_{T_k}$ as
			$$S_{T_k}:=\min_{major~i}\{S_{T}(X^k_{i}|C(X^k_{i})=C_{k}\}$$ Also, update the set $N_k $ of non-dominated pairs (till now) from which the original pair $(S_{T_k}, C_k)$ was selected and repeat the general step for the updated pair  $(S_{T_k}, C_k)$.\end{enumerate} }
	\item Construct a set $N_{k+1}$ consisting of non-dominated pairs (NDPs) till now obtained from the OMFSs of the  major sub-problems $P_{[a^k_i,b^k_i]}^i(S_{T_{k}})$ and  those obtained at previous steps (if any), i.e.,\vspace{2mm}\\ $N_{k+1}=\left\{\mbox{ NDPs obtained from major sub-problems} ~P_{[a^k_i,b^k_i]}^i(S_{T_{k}}) \right \}\cup$ \\$ \left\{\mbox{ NDPs in the set } N_k ~\mbox{if there is any} \right \} $\vspace{2mm}\\
	\noindent The repeated pairs are ignored while constructing the set $N_{k+1}$.
	Now, select a pair $(S_{T_{k+1}},C_{k+1})$ from the set $N_{k+1}$ where $C_{k+1}=\displaystyle\min_{(S_T, C)\in N_{k+1}}\{C\}$ and repeat the general step by setting $k=k+1$.\end{enumerate}
\textbf{Step 5.}
The non-dominated pairs are $(S_{T_0},C_0)$,$(S_{T_1},C_1)$,$\cdots$,$(S_{T_{k}},C_{k})$. \vspace{3mm}\\
The flow chart of the algorithm is given in Figure \ref{ttt}.
\begin{remark}
	While solving the problem $CP_0$, if one of the transportation  times $T^0_1$ or $T^0_2$ comes out to be zero, i.e., allocation is done only in one of the phases in order to get the optimal transportation cost, then this solution cannot be  considered as a feasible solution to the problem BTPTP (See mathematical formulation of   Phase-I and Phase-II problems). In case $T^0_1$ ($T^0_2$) takes zero value, then we solve the problem $CP_0$ again, with some links in Phase-II with transportation time greater than or equal to $T^0_2$ ($T^0_1$),  blocked. This is continued until the problem $CP_0$ provides   non-zero values of  $T^0_1$ and $T^0_2$.\\
	Further, if such a situation arises for a major sub-problem at any step of the algorithm, then this sub-problem is discarded and one of its immediate next minor sub-problem is considered as the major sub-problem for that particular step (if exists).
	Immediate next minor sub-problem of a major sub-problem is the one with time entries immediately lesser than that of the major sub-problem. Graphically, in Figure \ref{1}, the problem $P_{[6,6]}^{3.1}(14)$ is the immediate next minor sub-problem of the sub-problem $P_{[7,6]}^{3.2}(14)$.
	It should be noted that if we ignore the immediate next minor sub-problem at this step, then there may be a chance that this sub-problem does not occur in the  subsequent steps and hence, a non-dominated  pair may be missed while implementing the algorithm on an instance. \end{remark}
\begin{remark}Since, the sub-problems $P_{[a^k_i,b^k_i]}^i(S_{T_k})$ are constructed in such a way that  OMFSs of  these sub-problems (if exist)  yield a feasible pair corresponding to which the sum of transportation times is strictly less than $S_{T_k}$ and  the time entries in both the phases are assumed to be integers, therefore,  the sum of transportation times will be reduced  at least by one unit.\end{remark}
\begin{remark} At sub-step 3 of Step 3 of the BTPTP-Algorithm, we select a pair from the set $N_{k+1}$ which is non-dominated till now and has minimum transportation cost. However, if condition 2(a) of Step 3 holds for this pair, then it becomes the overall non-dominated  pair of the problem BTPTP. \end{remark}
\begin{remark} 
	The algorithm records a pair $(S_T,C)$ as a non-dominated pair if all the major sub-problem restricted at the sum $S_T$ that are M-feasible, yield transportation cost greater than $C$. The algorithm terminates as soon as it reaches a stage where all the major sub-problems corresponding to the pair $(S_{T_l},C_l)\in N_l$  where $N_{l}=\left\{\mbox{ NDPs obtained from major sub-problems} ~P_{[a^{l-1}_i,b^{l-1}_i]}^i(S_{T_{k}}) \right \}\cup$ \\$ \left\{\mbox{ NDPs in the set } N_{l-1} ~\mbox{if there is any} \right \} $, are non M-feasible. The pair $(S_{T_l},C_l)$ corresponds to the minimum cost pair among all the dominated pairs recorded at that particular instance. It means $S_{T_l}$ is the highest time among all the pairs in the set $N_l$ and $C_l$ is the minimum cost. If all the major sub-problems restricted at the sum $S_T$ are not M-feasible, then in this case, $N_l$ would be a singleton set containing $(S_{T_l},C_l)$ pair only. It is because if  $N_l$ is not singleton, then other pairs in the set would correspond to M-feasible solutions of some major sub-problems restricted at the sum $S_{T_l}$. This gives a contradiction as none of these major sub-problems is M-feasible. 	\end{remark}
\subsection{Theoretical development of the BTPTP-Algorithm}
\begin{lemma}\label{wahel} 
	Let an OMFS of a restricted CMTP solved at an iteration of the algorithm, yields a pair $(S_T,C)$ and  the optimal cost of one of the  sub-problems restricted at $S_T$  is equal to $C$. Then, the pair $(S_{T}, C)$ is not a non-dominated  pair.\end{lemma}
\noindent{Proof.} Let $(S_{T}, C)$ be a feasible pair obtained from an OMFS $X$, of a  restricted CMTP solved at any step of the algorithm i.e., $S_{T} = S_{T}(X)$ and $C =C(X)$.\\
At an OMFS ${{X}_{i}}$ of the sub-problem $P_{[a_i,b_i]}^{i}(S_{T})$,
\begin{equation}\label{2}
S_{T}({{{X}}_{i}}) < S_{T}(X),~ \forall~ \mbox{major values of} ~i.
\end{equation}
Now, if $CP_{[a_i,b_i]}^i(S_{T})=C$, at least for one major $i$,  then we can write 
\begin{equation}\label{33}
C({{{X}}_{i}})=C(X),  \mbox{~for~that~particular }~i.
\end{equation}
From the inequations $(\ref{2})$ and $(\ref{33})$, it follows that the pair $(S_{T}({{{X}}_{i}}), C({{{X}}_{i}}))$ dominates the pair $(S_{T}, C)$. Hence, $(S_{T}, C)$ is not a non-dominated  pair.
\begin{theorem}\label{wahet1}
If the optimal cost of all the M-feasible major sub-problems restricted at $S_T $ corresponding to  the pair $(S_T,C)$,  are greater than $C$, then the pair $(S_T,C)$ is a non-dominated pair. 
\end{theorem}

\begin{theorem}\label{wahet2}
Let the pair $(S_T,C)$ be such that none of the major sub-problems  restricted at $S_T$ is M-feasible. Then $(S_T,C)$ is  the last recorded efficient/non-dominated point with $S_T$ as the minimum sum of transportation times of both the phases.   
\end{theorem}
\begin{theorem}\label{wahet3}
The proposed algorithm    terminates in a finite number of steps and records all non-dominated  pairs i.e., no non-dominated  pair is missed by the algorithm.
\end{theorem}
Theorems \ref{wahet1},  \ref{wahet2} and \ref{wahet3} can be proved on the similar lines as done by  Kaur et al. \cite{Kaur} in their paper.
\section{Numerical Illustration}\label{Num}
Consider a $3 \times 5$ BTPTP 2 given in Table \ref{table1} where  cells with bold entries denote the Phase-I links while the others denote Phase-II links. Here, the index sets of sources and destinations are $I=\{1,2,3\}$ and $J=\{1,2,3,4,5\}$, respectively. Also, at each source-destination link, the number on top left corner denotes the transportation cost $c_{ij},~\mbox{for} ~i \in I,~ j\in J$ whereas the number on bottom right corner denotes the transportation time $t_{ij}$.  The quantities $a_i$ and $b_j$ denote the availability at $i^{th}$ source and demand at $j^{th}$ destination, respectively. We explain the implementation of the BTPTP-Algorithm to find all the non-dominated pairs of this problem.\vspace{2mm}\\
\begin{table}
\caption{The data of $3 \times 5$  BTPTP 2}
\centering
\begin{tabular}{c|c|c|c|c|c|c} \label{table1}
	& $D_1$ & $D_2$ & $D_3$ & $D_4$ & $D_5$ & $a_i$   \\
	\hline $S_1$ & \begin{tabular}{ll}
		$c_{11}=10$&~~\\& $t_{11}=5$
	\end{tabular}   & \begin{tabular}{lr}
		\bf\textbf{9} &\\ &\bf\textbf{6}
	\end{tabular}   &\begin{tabular}{lr}
		\bf\textbf{11} &\\ &\bf\textbf{3}
	\end{tabular}   &\begin{tabular}{lr}
		$7$ &\\ &$2$
	\end{tabular}   & \begin{tabular}{lr}
		$8$ &\\ &$3$
	\end{tabular}   &  7\\
	\hline
	$ S_2$   & \begin{tabular}{lr}
		\bf\textbf{11}~~~~ &~~\\~~ &~~~~\bf\textbf{2}
	\end{tabular}   &  \begin{tabular}{lr}
		$10$ &\\ &$3$
	\end{tabular}     &  \begin{tabular}{lr}
		\bf\textbf{13} &\\ &\bf\textbf{5}
	\end{tabular}   & \begin{tabular}{lr}
		$14$ &\\ &$8$
	\end{tabular}   &  \begin{tabular}{lr}
		\bf\textbf{12} &\\ &\bf\textbf{4}
	\end{tabular}    & 8   \\
	\hline
	$ S_3$ & \begin{tabular}{lr}
		$8$~~~~~~~ &\\~~~~~~~~~~~ &$4$
	\end{tabular}     & \begin{tabular}{lr}
		\bf\textbf{6} &\\ &\bf\textbf{5}
	\end{tabular}    & \begin{tabular}{lr}
		$9$ &\\ &$8$
	\end{tabular}   & \begin{tabular}{lr}
		\bf\textbf{10} &\\ &\bf\textbf{3}
	\end{tabular}   & \begin{tabular}{lr}
		$13$ &\\ &$8$
	\end{tabular}    &9\\
	\hline $b_j$ & 3 & 5 & 6 & 6 & 4 &
\end{tabular}
\end{table}
The time entries of Phase-I and Phase-II links can be arranged in ascending order as follows:\\Phase-I entries:
$~~~~~~t_1^1 (=2)<t_1^2 (=3) < t_1^3 (=4) < t_1^4 (=5) <t_1^5 (=6) $.\vspace{2mm}\\
Phase-II entries:
$~~~~~~t_2^1 (=2)<t_2^2 (=3) < t_2^3 (=4) < t_2^4 (=5) <t_2^5 (=8) $.\\

\begin{figure}[h]
\centering
\begin{tabular}{c}
	\begin{minipage}{5.0in}
		\centering
		\includegraphics[width=\textwidth]{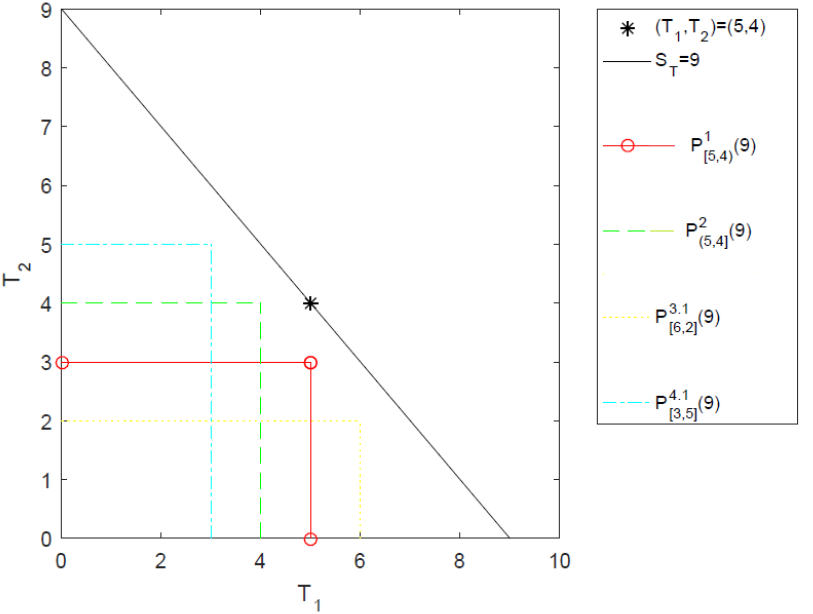}
	\end{minipage}
\end{tabular}
{\footnotesize\caption{Upper bounds for the time entries of  sub-problems restricted at the sum $9$.}\label{5}}
\end{figure}

\textbf{Step 1.}
\begin{enumerate}
\item Enter the number of sources, $m=3$.
\item Enter the number of destinations, $n=5$.
\item Generate the set of transportation links $I \times J$. 
\item Enter the availability of   sources as $7,8,9$.
\item Enter demand of destinations as $3,5,6,6,4$.
\item Construct the Phase-I and Phase-II sets namely, $PH_1=\{(1,2),(1,3),(2,1),(2,3),(2,5),(3,2),(3,4)\}$ and $PH_2=\{(1,1),(1,4),(1,5),(2,2),(2,4),(3,1),(3,3),(3,5)\}$.	
\end{enumerate}
\textbf{Step 2.}
\begin{enumerate}
\item Enter the time entries    $\{t_{ij}\}_{3 \times 5}$ as $\bigg[\begin{array}{ccccc}
	10	&2&5&9&8\\
	8&4&11&7&6\\
	3&5&4&6&9\\
\end{array}\bigg]$
\item Enter the cost entries   $C=\{c_{ij}\}_{3 \times 5}$ as $\bigg[\begin{array}{ccccc}
	3	&3&7&4&6\\
	2&6&4&9&9\\
	12&6&8&5&9\\
\end{array}\bigg]$	 
\end{enumerate}
\textbf{Step 3.}  Construct and solve the problem $CP_0$. An  OFS $X^{0}$ of this problem is obtained as \\
$x_{14}=6,~x_{15}=1,~x_{21}=3,~x_{23}=2,~x_{25}=3,~x_{32}=5,~x_{33}=4$. \\From $X^{0}$, the optimal cost $C_0$ is obtained as $211$ and the corresponding  feasible solutions $Y^0$ and $Z^0$, satisfying (\ref{15}), of Phase-I and Phase-II problems respectively, give $T_{1}(Y^0)=5$ and $T_{2}(Z^0)=8$. \\Therefore, the pair $(S_{T_0},C_0)$ is noted as $(13,211)$ where $S_{T_0}=T_{1}(Y^0)+T_{2}(Z^0)= 5+8=13$. Define $N_0= \{(S_{T_0},C_0)\}=\{(13,211)\}$. \vspace{1mm}\\
\begin{figure}[h]
\centering
\begin{tabular}{c}
	\begin{minipage}{5.0in}
		\centering
		\includegraphics[width=\textwidth]{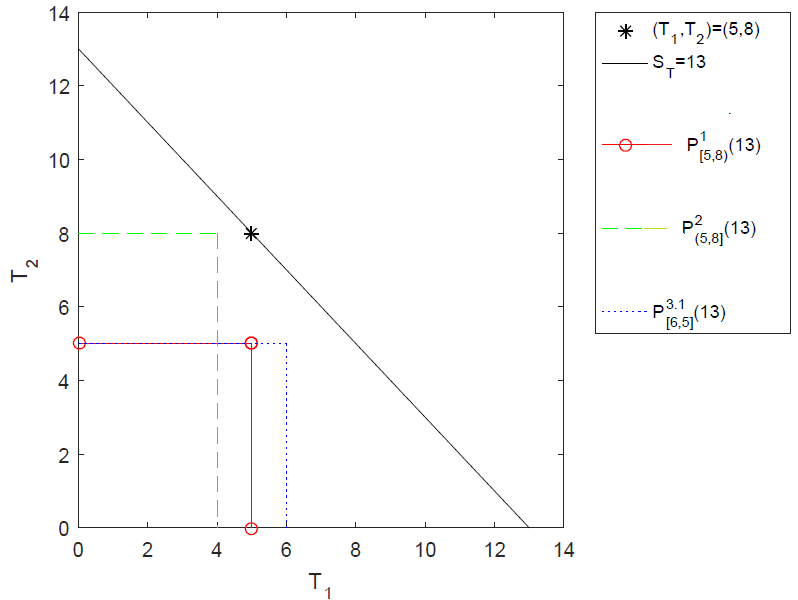}
	\end{minipage}
\end{tabular}
{\footnotesize\caption{Upper bounds for the time entries of  sub-problems restricted at the sum $13$. }\label{3}}
\end{figure}
\begin{figure}[h]
\centering
\begin{tabular}{c}
	\begin{minipage}{5.0in}
		\centering
		\includegraphics[width=\textwidth]{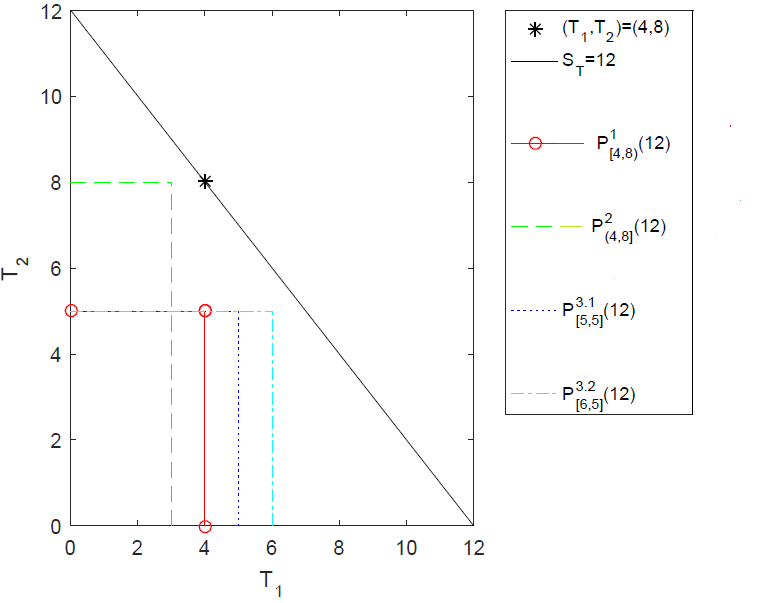}
	\end{minipage}
\end{tabular}
{\footnotesize\caption{Upper bounds for the time entries of various sub-problems restricted at the sum $12$.}\label{4}}
\end{figure}\\
\textbf{Step 4}.\\\textbf{Sub-step 1}
Define the sub-problems restricted at the sum $13$ i.e.,  $P_{[a^0_i,b^0_i]}^i(13)$,  $\forall~ i$ as  given below:  \begin{itemize}
\item{ $P_{[5,5]}^1(13) $ }
\item{ $P_{[4,8]}^2(13) $ }
\item{Here $T'_1=\displaystyle\min\{S_{T_{0}}, (T_{1})_{max}\}=\displaystyle\min\{13, 6\}=6$
	and there is no link in Phase-I set with the time entry  lying between $T_1(Y^0)=5$ and $T'_1=6$, i.e., $r=0$.  Therefore, $T_{1.1}$ and  $T'_{1}$ coincide. So, there is only one sub-problem viz.,
	$P_{[6,5]}^{3.1}(13)$.}
\item{Also, $T'_2=\displaystyle\min\{S_{T_{0}}, (T_{2})_{max}\}=\displaystyle\min\{13, 8\}=8$
	and there is no link in Phase-II set with the time entry  lying between $T_2(Z^0)=8$ and $T'_2=8$, i.e., $s=0$. Therefore, $T_{2.1}$ and  $T'_{2}$ coincide. So there is only one sub-problem viz., $P_{[4,8]}^{4.1}(13)$. }
\end{itemize} The  upper bounds for the time entries for the links of Phase-I and Phase-II sets corresponding to these sub-problems are depicted in the form of rectangles in $(T_1, T_2)$ plane, in Figure \ref{3} from which we observe that only the sub-problems $P_{[4,8]}^2(13) $ and  $P_{[6,5]}^{3.1}(13)$ are major sub-problems. The OMFSs of these sub-problems   are obtained as follows:\vspace{2mm}\\
$P_{[4,8]}^2(13) $: The OMFS $X^0_2$, of this sub-problem is obtained as\\
$x_{14}=6,~x_{15}=1,~x_{21}=0,~x_{22}=5,~x_{25}=3,~x_{31}=3,~x_{33}=6$
which gives $T_1(Y^0_2)=4$, $T_2(Z^0_2)=8$ and $C(X^0_2)= 214$. In other words, we say $CP_{[4,8]}^{2}(13)=214$.\\ Therefore, $(S_T(X^0_2), C(X^0_2))= (4+8, 214)=(12, 214)$. \vspace{2mm} \\
$P_{[6,5]}^{3.1}(13)$: The OMFS $X^0_{3.1}$, of this sub-problem is obtained as\\
$x_{14}=5,~x_{15}=2,~x_{23}=6,~x_{25}=2,~x_{31}=3,~x_{32}=5,~x_{34}=1$ which gives $T_1(Y^0_{3.1})=5$, $T_2(Z^0_{3.1})=4$ and $C(X^0_{3.1})= 217$. In other words, we say $CP_{[6,5]}^{3.1}(13)=217$.\\ Therefore, $(S_T(X^0_{3.1}), C(X^0_{3.1}))=(5+4, 217)=(9, 217)$.\vspace{2mm}\\
Thus, we see that $CP_{[a^0_i,b^0_i]}^{i}(13)> 211, ~\forall ~\mbox{major values of}~i$. Therefore, $(13,211)$ is the first (overall) non-dominated  pair.\\
\noindent Now, construct the set
$$\begin{array}{rl}N_1=&\left\{\mbox{ NDPs obtained by solving the major sub-problems} ~P_{[a^0_i,b^0_i]}^i(S_{T_{0}}) \right \} \cup \\& \left\{\mbox{NDPs obtained in previous steps} \right \}\\=&\{(12,214),(9,217)\}\cup\emptyset\end{array}$$
Out of the set $N_{1}$, select the pair $(12,214)$ as it has minimum sum of transportation costs. Call this pair as $(S_{T_{1}},C_{1})$ and construct the sub-problems $P_{[a^1_i,b^1_i]}^i(S_{T_{1}})$,  $\forall~ i$. Note that the pair $(S_{T_{1}},C_{1})$ is obtained from an OMFS $X^0_2$ of a sub-problem solved at the Step 1.    \\
\noindent\textbf{Sub-step 2}. Define the sub-problems restricted at the sum $12$ i.e., $P_{[a^1_i,b^1_i]}^i(12) $,  $\forall~ i$ as  given below:
\begin{itemize}
\item{ $P_{[4,5]}^1(12) $}
\item{ $P_{[3,8]}^2(12) $}
\item{  $P_{[5,5]}^{3.1}(12)$ and
	$P_{[6,5]}^{3.2}(12)$
}
\item{
	$P_{[3,8]}^{4.1}(12)$  }
\end{itemize}
The  upper bounds for the time entries of the  links of Phase-I and Phase-II sets corresponding to these sub-problems  are depicted in the form of rectangles in $(T_1,T_2)$ plane, in Figure \ref{4} from which we observe that the sub-problems   $P_{[3,8]}^2(12) $ and  $P_{[6,5]}^{3.2}(12)$ are major sub-problems. The OMFSs of these sub-problems are  obtained as follows:\vspace{1mm}\\
$P_{[3,8]}^2(12) $: The OMFS $X^1_2$ of this sub-problem is obtained as \\
$x_{14}=3,~x_{15}=4,~x_{21}=3,~x_{22}=5,~x_{24}=0,~x_{33}=6,~x_{34}=3$ which gives $T_1(Y^1_2)=3$, $T_2(Z^1_2)=8$ and $C(X^1_2)= 220=CP_{[3,8]}^{2}(12)$.\\ Therefore, $(S_T(X^1_2), C(X^1_2))= (3+8, 220)=(11, 220)$. This pair is dominated by the pair $(9,217)$.
\vspace{1mm}\\  $P_{[6,5]}^{3.2}(12)$: The  upper bounds for the time entries of the  links of Phase-I and Phase-II sets corresponding to this sub-problem is same as those of the sub-problem  $P_{[6,5]}^{3.1}(13)$ whose OMFS provided the pair $(9, 217)$ in Step 1, therefore, an OMFS of the sub-problem  $P_{[6,5]}^{3.2}(12)$ will also yield the same pair i.e., $(S_T(X^1_{3.2}), C(X^1_{3.2}))= (5+4, 217)=(9, 217)$. \vspace{1mm}\\
Thus, we see that $CP_{[a^1_i,b^1_i]}^{i} (12)> 214 ~\forall ~\mbox{major values of}~i$. Therefore, $(12,214)$ is the second (overall) non-dominated  pair.\\ Now, the set $N_{2}$ of non-dominated pairs (till now)  becomes
$$N_2=\emptyset\cup\{(9,217)\}$$
Now, there is only one pair in the set $N_2$, therefore, we take $(S_{T_{2}},C_{2})=(9,217)$ and  and construct the sub-problems $P_{[a^2_i,b^2_i]}^i(S_{T_{2}})$,  $\forall~ i$ .\vspace{1mm}\\
\textbf{Sub-step 3}.  Define the sub-problems restricted at the sum $9$, i.e., $P_{[a^2_i,b^2_i]}^i(9) $,  $\forall~ i$ as  given below.
\begin{itemize}
\item{ $P_{[5,3]}^1(9) $ }
\item{ $P_{[4,4]}^2(9) $  }
\item{
	$P_{[6,2]}^{3.1}(9)$  }

\item{
	$P_{[3,5]}^{4.1}(9)$ and $P_{(1,8]}^{4.2}(9)$ (This sub-problem is not possible to construct as there is no link in  Phase-I set with the time entry strictly less than 1)
	
}
\end{itemize}
The  upper bounds for the time entries of  links of Phase-I and Phase-II sets corresponding to  the above sub-problems (which are possible to construct)  are depicted in the form of rectangles in $(T_1,T_2)$ plane, in Figure \ref{5} from which we observe that the  sub-problems   $P_{[5,3]}^1(9) $, $P_{[4,4]}^2(9) $,
$P_{[6,2]}^{3.1}(9)$ and
$P_{[3,5]}^{4.1}(9)$ are major sub-problems. The OMFSs of these sub-problems are  obtained  as follows:\vspace{1mm}\\
$P_{[5,3]}^1(9) $: The OMFS $X^2_1$ of this sub-problem is \\
$x_{13}=1,~x_{14}=2,~x_{15}=4,~x_{21}=3,~x_{23}=5,~x_{32}=5,~x_{34}=4$ which gives $T_1(Y^2_1)=5$, $T_2(Z^2_1)=3$ and $C(X^2_1)= 225=CP_{[5,3]}^{1}(9)$.\\ Therefore, $(S_T(X^2_1), C(X^2_1))= (5+3, 225)=(8, 225)$.
\vspace{1mm}\\
$P_{[4,4]}^2(9) $: The OMFS $X^2_2$ of this sub-problem is\\
$x_{13}=6,~x_{14}=0,~x_{15}=1,~x_{22}=5,~x_{25}=3,~x_{31}=3,~x_{34}=6$ which gives $T_1(Y^2_2)=4$, $T_2(Z^2_2)=4$ and $C(X^2_2)= 244=CP_{[4,4]}^{2}(9)$.\\ Therefore, $(S_T(X^2_2), C(X^2_2))= (4+4, 244)=(8, 244)$.  This pair is dominated by the pair $(8,225)$. \\$P_{[6,2]}^{3.1}(9)$: The OMFS $X^2_{3.1}$ of this sub-problem is \\
$x_{13}=5,~x_{14}=2,~x_{21}=3,~x_{23}=1,~x_{25}=4,~x_{32}=5,~x_{34}=4$ which gives $T_1(Y^2_{3.1})=5$, $T_2(Z^2_{3.1})=2$ and $C(X^2_{3.1})=233=CP_{[6,2]}^{3.1}(9)$.\\ Therefore, $(S_T(X^2_{3.1}), C(X^2_{3.1}))=(5+2, 233)=(7, 233)$.\vspace{1mm}\\
$P_{[3,5]}^{4.1}(9)$: There is no M-feasible solution to this sub-problem.\vspace{1mm}\\
Thus, we see that $CP_{[a^2_i,b^2_i]}^{i} (9)> 217, ~\forall ~\mbox{major values of}~i$. Therefore $(9,217)$ is the third non-dominated  pair.
\begin{figure}[h]
\centering
\begin{tabular}{c}
	\begin{minipage}{5.0in}
		\centering
		\includegraphics[width=\textwidth]{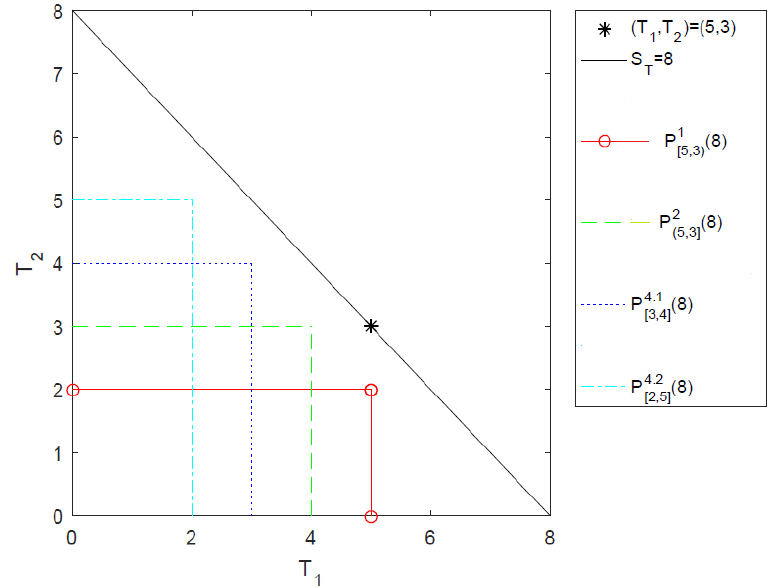}
	\end{minipage}
\end{tabular}
{\footnotesize\caption{Upper bounds for the time entries of  sub-problems restricted at the sum $8$.}\label{6}}
\end{figure}
Now, the set $N_{3}$ of non-dominated pairs becomes
$$N_3=\{(8,225)(7,233)\}\cup\emptyset$$
Select the pair $(8,225)$ as  it has the minimum sum of transportation costs and set $(S_{T_{3}},C_{3})=(8,225)$ and construct the sub-problems $P_{[a^3_i,b^3_i]}^i(8) $,  $\forall~ i$.\vspace{1mm}\\
\textbf{Sub-step 4}. Define the sub-problems restricted at the sum $8$ i.e., $P_{[a^3_i,b^3_i]}^i(8) $,  $\forall~ i$ as  given below.
\begin{itemize}
\item{ $P_{[5,2]}^1(8)$
}
\item{ $P_{[4,3]}^2(8)$}
\item{
	$P_{[6,2)}^{3.1}(8)$ (This sub-problem is not possible to construct as there is no link in  Phase-II set with the time entry strictly less than 2).}
\item{$P_{[3,4]}^{4.1}(8)$ and  $P_{[2,5]}^{4.2}(8)$ }
\end{itemize}
\begin{figure}[h]
\centering
\begin{tabular}{c}
	\begin{minipage}{5.0in}
		\centering
		\includegraphics[width=\textwidth]{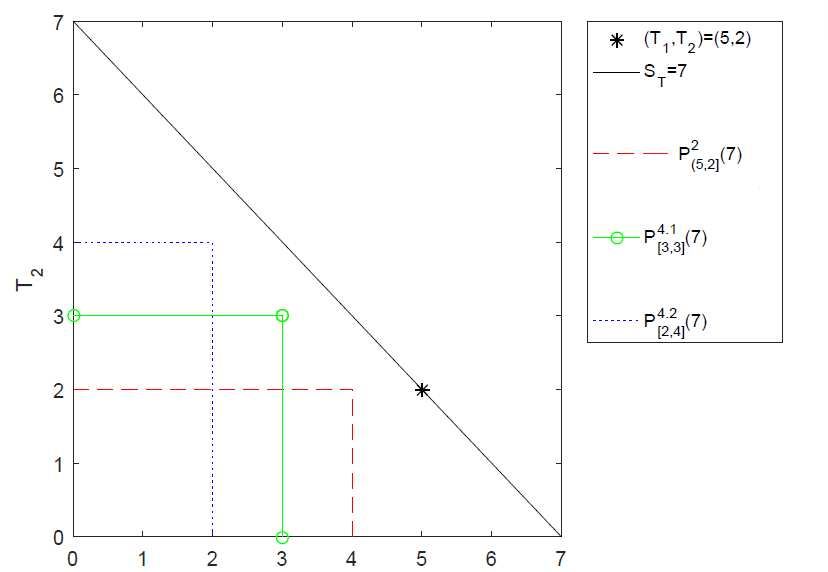}
	\end{minipage}
\end{tabular}
{\footnotesize\caption{Upper bounds for the time entries of  sub-problems restricted at the sum $7$.}\label{7}}
\end{figure}
The  upper bounds for the time entries of the links of Phase-I and Phase-II sets corresponding to  the above sub-problems (which are possible to construct)  are depicted in the form of rectangles in $(T_1,T_2)$ plane, in Figure \ref{6} from which we observe that the  sub-problems
$P_{[5,2]}^1(8)$, $P_{[4,3]}^2(8)$,
$P_{[3,4]}^{4.1}(8)$ and  $P_{[2,5]}^{4.2}(8)$ are major sub-problems. The OMFSs of these sub-problems  are obtained  as follows:\vspace{1mm}\\
$P_{[5,2]}^1(8)$: The OMFS $X^3_1$ of this sub-problem is obtained as  \\
$x_{13}=5,~x_{14}=2,~x_{21}=3,~x_{23}=1,~x_{25}=4,~x_{32}=5,~x_{34}=4$ which gives $T_1(Y^3_1)=5$, $T_2(Z^3_1)=2$ and $C(X^3_1)= 233=CP_{[5,2]}^{1}(8)$.\\ Therefore, $(S_T(X^3_1), C(X^3_1))= (5+2, 233)=(7, 233)$ (Repeated). \vspace{1mm}
\\
$P_{[4,3]}^2(8)$: There is no M-feasible solution to this sub-problem.\vspace{1mm}\\
The sub-problems $P_{[3,4]}^{4.1}(8)$ and $P_{[2,5]}^{4.2}(8)$ are contained in the sub-problem $P_{[3,5]}^{4.1}(9)$ which has no M-feasible solution. Therefore, these problems also, do not have an M-feasible solution.
Note that $CP_{[a^3_i,b^3_i]}^{i} (8)> 225, ~\forall ~\mbox{major values of}~i$ for which an OMFS exists. Therefore, the pair  $(8,225)$ is the fourth non-dominated  pair. Now, the set $N_{4}$ of non-dominated pairs (till now) becomes
$$N_4=\emptyset\cup\{(7,233)\}$$
Out of $N_{4}$, only choice is the pair $(7,233)$. Call it $(S_{T_{4}},C_{4})$ and construct the sub-problems $P_{[a^4_i,b^4_i]}^i(7)$,  $\forall~ i$. Note that the pair $(S_{T_{4}},C_{4})$ is initially obtained from the solution $X^2_{3.1}$ of a sub-problem solved at the step 3.\\
\textbf{Sub-step 5}. Define the sub-problems restricted at the sum $7$ i.e., $P_{[a^4_i,b^4_i]}^i(7)$,  $\forall~ i$ as given below.
\begin{itemize}
\item{ $P_{[5,2)}^1(7)$ (This sub-problem is not possible to construct as there is no link in  Phase-II set with the time entry strictly less than 2).}
\item{ $P_{(5,2]}^2(7)$}
\item{
	$P_{[6,1)}^{3.1}(7)$(This sub-problem is not possible to construct as there is no link in  Phase-II set with the time entry strictly less than 1). }
\item{
	$P_{[3,3]}^{4.1}(7)$,
	$P_{[2,4]}^{4.2}(7)$ and $P_{(2,5]}^{4.3}(7)$ (This sub-problem is not possible to construct as there is no link in  Phase-I set with the time entry strictly less than 2).}\end{itemize}
The  upper bounds for the time entries of  the links of Phase-I and Phase-II sets corresponding to  the above sub-problems (which are possible to construct)  are depicted in the form of rectangles in $(T_1,T_2)$ plane, in Figure \ref{7} from which we observe that the  sub-problems
$P_{(5,2]}^2(7)$, $P_{[3,3]}^{4.1}(7)$ and $P_{[2,4]}^{4.2}(7)$ are major sub-problems. None of these sub-problems is M-feasible. Therefore, there does not exist a  non-dominated pair with sum of transportation times strictly less than $7$. In other words,  $7= \displaystyle\min_{X \in S} \{S_T(X)\}$. \\
\textbf{Step 5. }  The BTPTP-Algorithm terminates as a non-dominated  pair with the minimum sum of transportation times of Phase-I and Phase-II has been  obtained.
Hence, the set of non-dominated  pairs is
$$\{(13,211), (12,214), (9,217),(8,225),(7,233)\}$$
The feasible solutions from which these non-dominated  pairs are obtained are considered as the Pareto optimal solutions of the $3 \times 5$ BTPTP 2.
All of these solutions are equally good as no additional subjective preference is considered in the construction of the problem.  \vspace{3mm}\\
\textit{Worst case complexity of the BTPTP-Algorithm:}
The problem BTPTP has been studied through various CMTPs which are defined according to the  upper bounds for the time entries of the links in Phase-I and Phase-II sets, at a particular step of the algorithm. These upper bounds  are obtained by exploring all the distinct  time entries in Phase-I and Phase-II such that the sum of times of both the phases corresponding to the optimal solutions of these problems, is strictly less than that obtained at the previous step. Therefore, the maximum number of CMTPs solved by the BTPTP-Algorithm to yield all the non-dominated points of the problem BTPTP is equal to the maximum possible combinations of Phase-I and Phase-II time entries. If the number of distinct time entries in Phase-I and Phase-II are $s_1$ and $s_2$ say, where $s_1+s_2\leq mn$,  then their maximum possible combinations are $s_1s_2$. Equivalently, the algorithm needs to solve at the most $s_1s_2$ CMTPs. Furthermore, the value of $s_1s_2$ for the given values of $m$ and $n$ is maximum when  $s_1=\Bigl\lfloor\frac{mn}{2}\Bigr\rfloor $ and $s_2=\left\lceil\frac{mn}{2}\right\rceil
$ or vice-versa where $\Bigl\lfloor.\Bigr\rfloor $ and $\Bigl\lceil.\Bigr\rceil $ represent the floor and ceiling functions, respectively.
Since the best time complexity to solve a CMTP of size $m \times n$, obtained by  Orlin \cite{orlin} is $O(\overline{m+n}+mn \log{(mn)})$   where $m$ and $n$ are the number of sources and destinations respectively, therefore the time complexity of the BTPTP-Algorithm is $O\Big(\Bigl\lfloor\frac{mn}{2}\Bigr\rfloor\left\lceil\frac{mn}{2}\right\rceil(\overline{m+n}+mn \log{(mn)})\Big)$.
\section{Computational Experiment}	
The algorithm given in Subsection \ref{algorithm} has been coded in MATLAB and experiments have been  conducted using Intel Processor i5 with 2.5 GHz, 4 GB RAM on 64-bit windows operating system. The algorithm  runs successfully for randomly generated problems of different sizes. The integral values of $t_{ij}$ and $c_{ij}$ in these problems have been drawn from a uniform distribution of the intervals $[1,40]$ and $[5,60]$ respectively. Table \ref{table2} narrates the computational behavior  of these problems with different percentages of Phase-I links. In a problem of  size $(m \times n)$, we report the average running time (in seconds) (taken over   50 instances) of the algorithm with different percentage of Phase-I links where the number of Phase-I links is taken as greatest integer value (non-zero) of the number corresponding to a particular percentage of $mn$. The running time of the algorithm depends upon the variation in time entries of Phase-I and Phase-II and the number of sub-problems solved for a particular instance of BTPTP, which further depends upon the sum of times corresponding to the minimum cost and the minimum possible sum of times. This is due to the reason that the algorithm is constructed in such a way that sum of times is reduced at each step, strictly by at least one unit. Therefore, greater the difference between sum of times corresponding to the minimum cost and the minimum possible sum of times, larger number of sub-problems are expected to be constructed which results in higher running time of the algorithm.
\par To analyze the computational behavior of the BTPTP-Algorithm more precisely, with respect to variations in time and cost entries 	
various   computational experiments have been conducted, in particular, for the $10 \times 10$ instance, by changing the interval from which the integral values of ${t_{ij}}$ and ${c_{ij}}$ have been drawn as follows:
\begin{enumerate}
\item{$t_{ij}\in [10,20], c_{ij}\in [5,60]$}
\item{$t_{ij}\in [10,20], c_{ij}\in [5,20]$}
\item{$t_{ij}\in [~1,40],  c_{ij}\in [5,20]$}
\end{enumerate}
Table \ref{table3} narrates the computational behavior of the algorithm for the $10 \times 10$
instance with the aforementioned variations from which we observe that if the variation in integral values of ${t_{ij}}$ and ${c_{ij}}$ is less, the running time of the algorithm is reduced significantly in comparison to that listed in Table \ref{table2}. For instance, a $10 \times 10 $ BTPTP  with $40 \%$ Phase-I links is solvable in $0.696354~secs$ with $t_{ij}\in [10,20], c_{ij}\in [5,60]$, which is almost $18\%$ of the running time of the algorithm for the same sized problem with  $t_{ij}\in [1,40], c_{ij}\in [5,60]$ (given in Table \ref{table2}). Thus, for the BTPTP instances with lesser variation in integral values of ${t_{ij}}$ and ${c_{ij}}$, our algorithm performs even better than it does for the instances with larger variations.
\begin{landscape}
\footnotesize{
	\begin{center}
		\begin{longtable}{ccccccccccc}
			\caption{Running time of BTPTP-Algorithm\label{table2}}\\
			\hline
			&      \multicolumn{9}{c}{Percentage of Phase-I links}\\ [0.5ex]
			\raisebox{1ex}{Size($m \times n$)}    & 10 &  20  &  30  &  40  &  50  &  60  &  70  &  80  &  90 \\ [0.5ex]
			\hline
			\endfirsthead
			\multicolumn{9}{c}%
			{\tablename\ \thetable\ -- {Continued from previous page}} \\
			\hline
			&     \multicolumn{9}{c}{Percentage of Phase-I links}\\ [0.5ex]
			\raisebox{1ex}{Size($m \times n$)}   & 10 &  20  &  30  &  40  &  50  &  60  &  70  &  80  &  90  \\ [0.5ex]
			\hline
			\endhead
			\hline \multicolumn{11}{c}{\textit{Running time of BTPTP-Algorithm}} \\
			\endfoot
			\hline
			\endlastfoot
			{3 $\times$ 3}     & 0.008
			& 0.010 &  0.017  & 0.01 & 	0.015 &	0.013  &	0.01	 &  0.006 &-     \\[0.5ex] 
			{3 $\times$ 4}
			&0.016	&0.024	&0.027 &	0.028 &	0.021 &	0.023 &	0.017 &	0.014 &	0.008  \\[0.5ex]

			{3 $\times$ 5}   & 0.02 	&0.033 	&0.025 	&0.031 &	0.029 &	0.030&	0.026 	&0.022 &	0.009   \\[0.5ex]
			{4 $\times$ 4 }     & 0.024
			& 0.05 &  0.033  & 	0.031  & 	0.038   & 0.034  & 	0.033  & 	0.025  & 	0.011
			\\[0.5ex]

			{3 $\times$ 6 }     & 0.029 &	0.052 &	0.044 &	0.046 &	0.046 &	0.044 &	0.048 	&0.035 &	0.012
			\\[0.5ex]
			{4 $\times$ 5 }     &0.029 &	0.057 &	0.068 	&0.066 &	0.035 &	0.037 &	0.064 	&0.042 	&0.027
			\\[0.5ex]
			{4 $\times$ 6 }     & 0.06 &	0.120 	&0.037 	&0.056 	&0.057 &	0.049 &	0.050 &	0.034 &	0.021
			
			\\[0.5ex]
			{5 $\times$ 5 }     & 0.086 &	0.040&		0.068	&	0.139&		0.104&		0.121	&	0.099&		0.091&		 0.042
			\\[0.5ex]
			{4 $\times$ 7 }     & 0.059
			& 0.096 &	0.109  &	0.173  &	0.206  &	0.144  &	0.103  &	0.107 &	0.079   \\[0.5ex]
			{5 $\times$ 6}     & 0.151
			& 0.117  &	0.112 &	0.172  &	0.274 & 	0.183  &	0.112 &  0.081  &		0.095    \\[0.5ex]
			{4 $\times$ 8}     & 0.179
			& 0.117  &	0.112 &	0.172  &	0.274 & 	0.183  &	0.112 &  0.081  &		0.095    \\[0.5ex]
			{5 $\times$ 7}     &0.388	&0.328	&0.455&	0.445&	0.465	&0.405&	0.353&	0.236&	0.094
			\\[0.5ex]
			{6 $\times$ 6}     & 0.124&	0.599&	0.509&	0.414&	0.502&	0.470&	0.365&	0.329	&0.123
			\\[0.5ex]
			{5 $\times$ 8}     & 0.429 &	0.748&	0.340 &	0.368 &	0.350 &	0.239 &	0.262	 &0.212 &	0.294
			\\[0.5ex]
			{6 $\times$ 7}     &0.448&	0.408	&0.315&	0.283&	0.253&	0.332&	0.231	&0.267&	0.146
			\\[0.75ex]
			{5 $\times$ 9}     & 0.753	&0.482&	0.325&	0.899&	0.436&	0.419&	0.260	&0.563	&0.155
			
			\\[0.75ex]
			{6 $\times$ 8}     & 1.003 	&1.628 &	0.608 	&0.725 &	0.594 &	0.796 &	0.641 &	0.293&	0.160
			\\[0.75ex]{5 $\times$ 10}     & 0.866 	&1.213 	&0.406 	&0.811 	&0.447 &	0.377 &	1.129 &	0.279 &	0.252
			\\[0.5ex]
			{8 $\times$ 7}     & 1.039 	&1.823 	&1.733 &	2.009 &	1.784 &	1.262 &	1.314 	&1.009 &	0.901
			\\[0.5ex]
			{10 $\times$ 6}     & 1.008 &	2.145 &	0.398 	&0.715 	&1.075 	&1.665 &	0.577 	&0.499 &	1.253
			\\[0.5ex]
			{8 $\times$ 8}     & 1.138 	&1.290 &	1.271 &	1.187 	&0.923 	&0.916 	&0.764 	&0.617 	&1.662
			\\[0.5ex]
			{8 $\times$ 9}     & 2.462    &  1.630 &	1.468 &	2.268 &	2.007 &	1.828 &	1.511 	&1.086 &	2.244     \\[0.5ex]
			
			{8 $\times$ 10}     & 1.656 &	2.059 &	2.954 &	3.187 &	2.605 &	1.960 &	1.592 &	1.501 &	2.577
			\\[0.5ex]
			{9 $\times$ 9}     & 1.383 	&6.767 &	1.88 &	2.177    &	2.193 &	1.848 &	1.497 	&1.172 &	2.280
			\\[0.45ex]
			{9 $\times$10}     & 3.527 	&2.198 &	2.989 &	11.091 &	11.189 	&7.924 &	7.693 	&2.335 	&1.419
			\\[0.45ex]
			{10 $\times$ 10}      &4.660 	&4.637 	&3.991 &	3.901 &	4.124 &	4.161 &	3.163 &	2.595 	&1.540
			
			\\[0.445ex]
			{20 $\times$ 20}     &95.152&	333.266&	239.424&	190.034	&183.847&	172.724&	157.747	&111.279&	115.197
			\\[0.45ex]
			{30 $\times$ 30}     &  518.731	&1079.501&	984.773&	1512.605	&2176.059&	1202.512	&1438.339&	1430.145	 &1207.578
			\\[0.45ex]
			{40 $\times$ 40}     &  8017.830	& 5396.779 & 3891.901	 &	4600.331 	& 4483.980 & 6607.785	 	&8457.318 &	 8196.052	& 5381.283
			\\[0.45ex]
			{50 $\times$ 50}     &   2.7567e+04	& 2.2835e+04 & 2.7553e+04	 &	 2.6589e+04	&  3.3617e+04 &	 2.1601e+04	& 2.5863e+04 &	2.4784e+04	& 1.5369e+04
			\\[0.45ex]
			\hline
		\end{longtable}
\end{center}}
\end{landscape}

\footnotesize{\begin{landscape}

	\begin{center}
		\begin{longtable}{cccc}
			\caption{Running time for 10 $\times$ 10 BTPTP for various time and cost distributions\label{table3}\label{baaaa}}\\
			\hline
			&      \multicolumn{3}{c}{Intervals}\vspace{3mm}\\
					Percentage of Phase-I links   & Time ~~~ Cost &  Time ~~~ Cost   & Time ~~~ Cost    \\ [0.5ex]
			& $[10,20]$~~$[5,60]$ &  $[10,20]$~~$[5,20]$  &  $[1,40]$~~$[5,20]$  \\ [0.5ex]\hline
			10  &0.528 &0.626 & 2.225  \\ [0.5ex]			
			20     &     0.463  &      0.536   &   2.907  \\ [0.5ex]
			30     &     0.718  &     0.794    &   4.33   \\ [0.5ex]
			40   &     0.828  &0.696 &2.540   \\ [0.5ex]
			50&0.828  &0.758 & 5.019 \\ [0.5ex]
			60&0.730  &0.721& 3.860  \\ [0.5ex]
			70 & 0.784&0.636&  2.126  \\ [0.5ex]
			80 &0.453  &0.665 & 2.335 \\ [0.5ex]
			90& 0.264&0.331 &1.270   \\ [0.5ex]
			\hline
		\end{longtable}
	\end{center}
\end{landscape}}

\section{Concluding Remarks}
\begin{enumerate}
\item  In this paper, an iterative algorithm is developed to find all the non-dominated pairs of a bi-objective two-phase transportation problem with the objective of minimizing `sum of transportation times' and `sum of transportation costs' of both the phases. The proposed algorithm at each iteration, solves various restricted CMTPs to find all the non-dominated solutions of the problem. 
\item The proposed algorithm terminates in a finite number of steps as the distinct time entries in Phase-I and Phase-II of the problem are finite.  
\item  The proposed algorithm does not miss any non-dominated solution  of the problem BTPTP and therefore,  solves the problem successfully without any limitation, however, a scope of improvement in the working of the algorithm may not be denied. For instance, the BTPTP-Algorithm  initially finds the optimal sum of transportation costs of both the phases and obtains the corresponding sum of transportation times and then reduces the latter   systematically, however, if we concentrate on finding the minimum sum of transportation times and the corresponding sum of transportation costs first  and then reducing the latter systematically, then it may be possible to observe better computational results for some BTPTP instances while the list of non-dominated points would be the same.			   
\item 	The computational behavior of the BTPTP-Algorithm in terms of CPU times,
depends on the variations in the time entries and the cost entries corresponding to a fixed size problem.  It has been observed from the Table \ref{baaaa} depicting this behavior that lesser the variation in the integral values of transportation times and costs, lesser the CPU time (for most of the instances) consumed by the algorithm.

\item It is important to note that the current problem has been discussed by assuming transportation times and transportation costs of all the source-destination links as crisp numbers, however, the practical scenario might be more general as due to many issues like bad weather, increased  fuel costs, traffic jams, strikes, waiting period at the destinations before unloading the goods etc., the  transportation times and  costs may vary significantly. In such  situations,  as a future aspect,  the current  problem   can be discussed  under uncertain environment and such a study would prove to be more connected with the real life problems. Another future aspect of this problem would be to discuss it as a multi-objective transportation problem, for instance, along with minimizing  total transportation time and cost, one might be willing  to optimize the deterioration cost of goods, quantity of goods delivered, under used capacity, reliability of delivery, energy consumption, safety of delivery and many other such objectives. Lastly, the study of the current problem can  be extended to a   multi-level optimization problem, solid transportation problem or as a transportation problem with non-linear costs.
\end{enumerate}

\end{document}